\documentclass{amsart} 
\usepackage{amsmath,amsthm,amssymb,tabularx}
\usepackage[cmtip,arrow,all,2cell]{xy}

\theoremstyle{plain}
\newtheorem{tm}{Theorem}[section]
\newtheorem{lm}[tm]{Lemma}
\newtheorem{cor}[tm]{Corollary}
\newtheorem{prop}[tm]{Proposition}

\theoremstyle{definition}
\newtheorem{definition}[tm]{Definition}
\newtheorem{ex}[tm]{Example}
\newtheorem{rem}[tm]{Remark}
\newtheorem{nots}[tm]{Notations}
\allowdisplaybreaks

\newcommand{\beq}{\begin{equation}}
\newcommand{\eeq}{\end{equation}}
\newcommand{\bga}{\begin{gather*}}
\newcommand{\ega}{\end{gather*}}
\newcommand{\bal}{\begin{align*}}
\newcommand{\eal}{\end{align*}}
\newcommand{\bit}{\begin{itemize}}
\newcommand{\eit}{\end{itemize}}
\newcommand{\btm}{\begin{tm}}
\newcommand{\etm}{\end{tm}}
\newcommand{\blm}{\begin{lm}}
\newcommand{\elm}{\end{lm}}
\newcommand{\bcor}{\begin{cor}}
\newcommand{\ecor}{\end{cor}}
\newcommand{\bprop}{\begin{prop}}
\newcommand{\eprop}{\end{prop}}
\newcommand{\bex}{\begin{ex}}
\newcommand{\eex}{\end{ex}}
\newcommand{\bpr}{\begin{proof}}
\newcommand{\epr}{\end{proof}}
\newcommand{\brem}{\begin{rem}}
\newcommand{\erem}{\end{rem}}
\newcommand{\bdf}{\begin{definition}}
\newcommand{\edf}{\end{definition}}
\newcommand{\bnots}{\begin{nots}}
\newcommand{\enots}{\end{nots}}

\def\C{\mathbb{C}}
\def\N{\mathbb{N}}
\def\R{\mathbb{R}}

\def\Id{\mathbb I}

\def\e{\varepsilon}
\let\a\alpha

\let\cal\mathcal
\def\m{\mathfrak m}

\def\id{{\rm id}}

\def \le {\leqslant}
\def \ge {\geqslant}
\let \<\langle
\let \>\rangle
\let\phi\varphi
\let\kappa\varkappa

\def\cst{$C^*$}
\def\Mhat{\widehat M}
\let\hat\widehat

\def\proj{\mathop{\widehat\otimes}}
\def\ohaag{\mathop{\overset{h}{\otimes}}}

\begin{document}
\author{Yulia Kuznetsova}
\title[A duality of locally compact groups]{A duality of locally compact groups that does not involve the Haar measure}
\address{University of Franche-Comt\'e, 16 route de Gray, 25030 Besan\c con, France}
\email{yulia.kuznetsova@univ-fcomte.fr}

\maketitle

{\obeylines\flushright
To the memory of my father,
Nikolay V.~Kuznetsov
}

\setcounter{tocdepth}{1}

\begin{abstract}
We present a simple and intuitive framework for duality of locally compacts groups, which is not based on the Haar measure. This is a map, functorial on a non-degenerate subcategory, on the category of coinvolutive Hopf \cst-algebras, and a similar map on the category of coinvolutive Hopf-von Neumann algebras. In the \cst-version, this functor sends $C_0(G)$ to $C^*(G)$ and vice versa, for every locally compact group $G$. As opposed to preceding approaches, there is an explicit description of commutative and co-commutative algebras in the range of this map (without assumption of being isomorphic to their bidual): these algebras have the form $C_0(G)$ or $C^*(G)$ respectively, where $G$ is a locally compact group. The von Neumann version of the functor puts into duality, in the group case, the enveloping von Neumann algebras of the algebras above: $C_0(G)^{**}$ and $C^*(G)^{**}$.
\end{abstract}

\section{Introduction}

A well-known theorem of L. S. Pontryagin states that a commutative locally compact group $G$ is isomorphic to its second dual group, where the dual group $\widehat G$ is the set of all unitary characters of $G$, which has a natural group structure.

This symmetry does not appear in the non-commutative case, because the set $\widehat G$ of irreducible unitary representations of $G$ --- the natural analogue of characters --- has no reasonable group structure. The problem can be stated, however, in an abstract form as follows.

Let $\cal H$ be a category. Call a functor $\hat{\ }:\cal H\to\cal H$ a duality if $\hat{\ \;}^{\,2}$ is equivalent to the identity functor. Let $\cal{LCG}$ and $\cal{LCAG}$ be the categories of all (respectively abelian) locally compact groups. On $\cal{LCAG}$ we have the Pontryagin duality functor $\hat{\;\;}$. The task is to construct a category $\cal H$ with duality, and a faithful functor $\cal A:\cal{LCG}\to \cal H$ so that the Pontryagin duality is preserved: $\widehat{\cal A(G)} \simeq \cal A(\widehat G)$ for every $G\in\cal{LCAG}$. This is illustrated by the following diagram:
$$
\xymatrix @R=1.pc @C=2.pc
 {
 \cal H\ar@{->}[rr]^{\hat{}}
 & &
\cal H
\\&&
 \\
 \cal{LCG} \ar@{->}[uu]^{\cal A}
 & &
 \cal{LCG}\ar@{->}[uu]^{\cal A}
\\&&
 \\
 \cal{LCAG} \ar@{->}[uu]^{}\ar@{->}[rr]^{\hat{}}
 & &
 \cal{LCAG}\ar@{->}[uu]^{}
 }
$$

In this setting, the duality problem has been solved by the theory of Kac algebras (the canonical reference is \cite{enock}), a theory later developed to that of locally compact quantum groups (see a recent survey monograph \cite{hajac}). A conceptual difference between this theory and that of Pontryagin is the distinguished role of the Haar measure (or its generalization called the Haar weight). In the classical theory, the dual group is defined in purely algebraic and topological terms: it is just the group of continuous characters. On the contrary, in the theory of Kac algebras the Haar weight is a part of the definition and is crucial in the construction of the dual object.

There has been a number of duality theorems which did not use the Haar measure explicitly. One alternative is the approach of multiplicative unitaries, developed mainly by S.~Baaj and G.~Skandalis \cite{baaj}, S.~L.~Woronowicz, T.~Masuda and Y.~Nakagami \cite{woron}, \cite{mn}, \cite{woron-mn}. These results have in common the fact that they are not constructive; formulating sufficient conditions for a duality, they do not present a means of obtaining reflexive objects (except for deriving them from a Haar weight).

Another alternative, very close to the present paper, is the work of E.~Kirchberg \cite{kirchberg} and its development by J.~Kustermans \cite{kustermans} and Ch.-K.~Ng \cite{ng}. The main idea is to define a dual algebra on the basis of ``unitary'' corepresentations. In the case when a Haar weight exists, one can show that the second dual algebra is isomorphic to the initial one. It would be very attracting to know that the first dual of any algebra is isomorphic to the third dual: this would give then a means of constructing reflexive algebras from any given ones. But for the moment, there are no such theorems.

In this paper, we propose another functor, close to the functor of Kirchberg, which has the following advantages. First, it is more explicit and accordingly more easy to calculate. Second, one can show that every commutative or cocommutative algebra in the range of our functor comes from a locally compact group and as a consequence is reflexive. This provides ground to a conjecture that the dual of every algebra is reflexive.

The main work is done in the category $\cal H_0$ of coinvolutive Hopf--von Neumann algebras (see definitions in section 2). In Section \ref{dual-von-neumann} we define a map $\hat{\;\;}$ on $\cal H_0$ and a full subcategory $\cal H$ of $\cal H_0$ on which this is a duality functor. The objects of $\cal H$ are just $M\in\cal H_0$ such that $M\simeq\hat{\hat M}$; let us call such algebras {\it reflexive}.

It is shown (Theorem \ref{von-neumann-duality}) that there is a faithful contravariant functor $\cal A: \cal {LCG}\to \cal H$, such that for an abelian group $G$, we have $\hat{\cal A(G)} \simeq \cal A(\hat G)$.
Explicitly, the functor $\cal A$ is given by $\cal A(G)=C_0(G)^{**}$. This is the enveloping von Neumann algebra of $C_0(G)$, canonically identified with the second dual space. The algebra $\hat{\cal A(G)}$ is the big group algebra of J.~Ernest, $W^*(G)\simeq C^*(G)^{**}$. The explicit consideration of the group case is given in Section \ref{reps-M(G)}.

Conversely, every commutative or co-commutative algebra in the range of our map is isomorphic to $C_0(G)^{**}$ or $C^*(G)^{**}$ respectively for some locally compact group $G$\, (Theorems \ref{commutative-is-group} and \ref{cocommutative-is-group}). All theorems of this kind known before assumed stronger properties guaranteeing a certain duality (a Haar weight, a regular representation etc).

In Section \ref{dual-cstar}, we give a \cst-algebraic version of this theory: for every coinvolutive \cst-bialgebra $A$ (see definition in Section \ref{dual-cstar}), one can define its dual coinvolutive \cst-bialgebra $\hat A$, such that $\hat{C_0(G)}=C^*(G)$ and $\hat{C^*(G)}=C_0(G)$. Every commutative or co-commutative algebra in the range of this map is isomorphic to $C_0(G)$ or $C^*(G)$ respectively.

\medskip

{\bf Acknowledgements}. I am grateful to Prof. M.~Enock for pointing me out to the thesis \cite{kirchberg} and for sending a copy of it. I thank Prof. E.~Kirchberg for clarification of some of his results. I am greatly obliged to the referee who pointed out several gaps and corrected numerous details in the paper.

\section{Definitions and notations}

\bnots\label{nots-tensor-products}
In general, we allow a von Neumann algebra to be zero, i.e. its unit may equal zero.
For a pair of von Neumann algebras $M,N$, we denote by $M\bar\otimes N$ their von Neumann tensor product. If $A$, $B$ are \cst-algebras, then $A\bar\otimes B$ denotes their spatial (minimal) tensor product. There should be no reason to confuse this case with the von Neumann algebras case. The dual of a Banach space $X$ is denoted by $X^*$. By $X\widehat\otimes_{\rm op}Y$, $X\ohaag Y$  we denote the (completed) projective operator space tensor product and the Haagerup tensor product of operator spaces $X,Y$ (more details on operator spaces can be found in \cite{effros}).
If $H$, $K$ are Hilbert spaces, then $H\otimes K$ is their Hilbert space tensor product, $B(H)$ is the space of bounded linear operators on $H$, and $N(H)$ is the space of trace class (nuclear) operators. Unless otherwise stated, $G$ will denote a locally compact group.
\enots

Recall the notion of the \cst-enveloping algebra, or the \cst-envelope of a Banach *-algebra \cite[\S11.1]{palmer}:

\bdf
Let $A$ be a Banach *-algebra, and let $I\subset A$ be the common kernel of all its *-representations (this is a two-sided *-ideal). $A$ is called {\it reduced\/} if $I=\{0\}$. For any $x\in A$, set $\|x\|_* = \sup\{ \|\pi(x)\|\}<\infty$, where supremum is taken over all *-representations $\pi$ of $A$. It is known that $\|x\|_*\le \|x\|$ for all $x\in A$. On $A/I$, the quotient seminorm is a norm; the completion of $A/I$ with respect to this norm is called the {\it \cst-envelope } of $A$ and denoted by $C^*(A)$. The canonical map from $A$ to $C^*(A)$ is injective if and only if $A$ is reduced.
\edf

From now on, the term ``representation'' will always mean ``*-representation''.
To every \cst-algebra $A$, one can associate in a canonical way a von Neumann algebra $W^*(A)$, which is called the {\it enveloping von Neumann algebra of $A$}  \cite[12.1.5]{dixmier} and has the following universality property:

\bprop
Let $A$ be a \cst-algebra and let $W$ be its von Neumann envelope. Let $\Phi:A\to W$ be the canonical morphism. Then for every representation $\pi: A\to B(H)$ there is unique normal representation $\tilde\pi: W\to B(H)$ such that $\tilde\pi(\Phi x)=\pi(x)$ for every $x\in A$. Moreover, $\tilde\pi(W)$ is the weak closure of $\pi(A)$. We will say that $\tilde\pi$ is the lifting of $\pi$ to $W$.
\eprop

Explicitly, $W^*(A)$ is constructed as the weak closure of the image of $A$ under the universal representation. Moreover, it can be naturally identified with the second dual space of $A$. Therefore we often write just $A^{**}$ instead of $W^*(A)$, having in mind this enveloping algebra structure on $A^{**}$.

A Banach *-algebra (strictly speaking, its image) is norm dense in its \cst-envelope, and a \cst-algebra is ultraweakly dense in its von Neumann envelope. We will also write $W^*(A)$ instead of $W^*(C^*(A))$ if $A$ is a Banach *-algebra.

The following definition can be found, e.g., in \cite[\S1.2]{enock}:

\bdf
A {\it coinvolutive Hopf-von Neumann algebra} is a triple $(M,\Delta,\kappa)$, where $M$ is a von Neumann algebra, $\Delta: M\to M\bar\otimes M$ (comultiplication) is an injective normal unital *-homomorphism such that $(\Delta\otimes\id)\Delta=(\id\otimes\Delta)\Delta$, and $\kappa:M\to M$ (coinvolution) is a *-antihomomorphism such that $\kappa^2=\id$ and $(\kappa\otimes\kappa)\Delta=\theta\Delta\kappa$, where $\theta$ is the flip map: $\theta(a\otimes b)=b\otimes a$.

A morphism of coinvolutive Hopf--von Neumann algebras $M,N$ is a normal *-homomorphism $\phi:M\to N$ such that: $\Delta_N\circ\phi=(\phi\otimes\phi)\Delta_M$ and $\kappa_N\circ\phi=\phi\circ\kappa_M$. Note that, contrary to \cite{enock}, we do not require that $\phi(1)=1$.
\edf

In the theory of operator spaces, there are two most natural notions of an algebra:

\bdf
A {\it completely contractive Banach algebra} is a Banach algebra $A$ which is an operator space such that the multiplication in $A$ is completely contractive, i.~e.~is extended to a continuous map $\m: A\proj_{\rm op}A\to A$. An {\it operator algebra} is a Banach algebra $A$ which is an operator space such that the multiplication in $A$ is extended to a continuous map $\m: A\ohaag A\to A$.
\edf

On the predual $M_*$ of a coinvolutive Hopf--von Neumann algebra $M$ one can introduce an involution, as usual in the Hopf theory: $\mu^*(a)=\overline{\mu(\kappa(a^*))}$, $\mu\in M_*$, $a\in M$. Then $M_*$ becomes a completely contractive Banach *-algebra, but in general not an operator algebra. For example, for the most popular algebra $M=L_\infty(G)$ one has $M_*=L_1(G)$, and this is known not to be an operator algebra.

When we speak of Banach *-algebras, we always suppose that the involution is isometric. It is known that the coinvolution $\kappa$ is also always isometric.

\subsection{Common group algebras}

There is a variety of algebras associated to a locally compact group $G$. We recall them here in order to have the freedom to use the notations below without extra explanations.

The most popular commutative algebras are: $C_0(G)$ --- the algebra of continuous functions vanishing at infinity; $C_b(G)$ --- the algebra of continuous bounded functions; $L_\infty(G)$ --- the algebra of equivalence classes of essentially bounded measurable functions; $A(G)$ --- the Fourier algebra, equal to the space of coefficients of the regular representation; $B(G)$ --- the Fourier--Stieltjes algebra, equal to the linear span of all continuous positive-definite functions. All these algebras are considered with pointwise multiplication, involution being the complex conjugation.

There is also a large choice of convolution algebras (and their completions): $M(G)$ --- the algebra of finite regular complex measures; $L_1(G)$ --- the subalgebra (in fact, an ideal) of absolutely continuous measures in $M(G)$; $C^*(G)$ --- the full group \cst-algebra, equal to the \cst-envelope of $L_1(G)$; $C^*_r(G)$ --- the reduced group \cst-algebra, generated by the regular representation of $L_1(G)$; $\cal L(G)$ --- the group von Neumann algebra, equal to the weak closure of $C^*_r(G)$ in $B(L_2(G))$; $W^*(G)$ --- the Ernest algebra, equal to the von Neumann envelope of $C^*(G)$.

In this list, $L_\infty(G)$, $\cal L(G)$ and $W^*(G)\simeq C^*(G)^{**}$ are well-known to have structures of coinvolutive Hopf--von Neumann algebras. There is also an algebra which is rarely used but is important in the sequel: $C_0(G)^{**}=W^*(C_0(G))$, the enveloping von Neumann algebra of $C_0(G)$. Since it can be identified with the second dual of $C_0(G)$, it is the dual space of $M(G)=C_0(G)^*$.

Instead of proving explicitly that $M(G)^*$ has a structure of a coinvolutive Hopf--von Neumann algebra, we can apply the known theory \cite[\S\,1.6]{enock} to $M=\cal L(G)$: its predual is $M_*=A(G)$, which has the \cst-envelope equal to $C_0(G)$ --- see a proof a few lines below --- and then $W^*(M_*)=C_0(G)^{**}$, as proved in \cite{enock}, has a structure of a coinvolutive Hopf--von Neumann algebra (agreeing with the structure of $C_0(G)$).

For future references it is convenient to formulate the following proposition, certainly known.

\bprop\label{envelope-A(G)}
$C^*(A(G))=C_0(G)$.
\eprop
\bpr
It is known that $A(G)$ is contained and dense in $C_0(G)$. The irreducible representations of $A(G)$ are just characters; it is known that every (nonzero) character of $A(G)$ has form $f\mapsto f(t)$, $f\in A(G)$, for some $t\in G$. It follows that $\|f\|_*=\sup_{t\in G} |f(t)|=\|f\|_\infty$ for every $f\in A(G)$. One can see that $\|\cdot\|_*$ is a norm (not just a seminorm), so that $C^*(A(G))$ is the completion of $A(G)$ with respect to $\|\cdot\|_*$. Now it is clear that $C^*(A(G))$ is just the closure of $A(G)$ in $C_0(G)$, and the statement follows.
\epr

\subsection{Multiplier algebras}\label{section-multipliers}

For a \cst-algebra $A$, let $M(A)$ denote the \cst-algebra of its two-sided multipliers \cite[\S~3.12]{pedersen}. It can be identified with a unital norm closed subalgebra in $W^*(A)$. A homomorphism $\phi: A\to M(B)$ is called non-degenerate if for an approximate identity $e_\a$ of $A$, $\phi(e_\a)$ converges to 1 in the strict topology of $M(B)$ (as operators on $B$). Every homomorphism $\phi:A\to M(B)$  has a unique extension to a $(A^*,B^*)$--weakly continuous homomorphism $\bar\phi: M(A)\to M(B)$. If $\phi$ is non-degenerate, $\bar\phi$ is unital. This applies to anti-homomorphisms as well.

For $A=C_0(G)$, $M(A)$ is the algebra $C_b(G)$ of bounded continuous functions on $G$. As a particular case of the discussion above, $C_b(G)$ is a *-subalgebra in $M(G)^*$, with the natural pairing $f(\mu) = \int f\,d\mu$, $f\in C_b(G)$, $\mu\in M(G)$.

\section{Representations with generator}

Let $M$ be a coinvolutive Hopf-von Neumann algebra, $M_*$ its predual. Fix a Hilbert space $H$ such that $M\subset B(H)$.  For $x,y\in H$ denote by $\mu_{xy}\in M_*$ the functional $\mu_{xy}(a)=\langle ax,y\rangle$, $a\in M$. Since  $M_*$ is a quotient space of $N(H)$, for every $\mu\in M_*$ there is a representation $\mu = \sum \mu_{x_n,y_n}$ with $x_n,y_n\in H$ such that $\sum \|x_n\|\,\|y_n\|<\infty$. Moreover, $\|\mu\|=\inf \sum\|x_n\|\,\|y_n\|$, where the infinum is taken over all such decompositions.

For another Hilbert space $K$, there is a natural isomorphism of operator spaces \cite[7.2.4]{effros} $(M_*\widehat \otimes_{op} N(K))^*\simeq M\bar\otimes B(K)$. In particular, every $U\in M\bar\otimes B(H)$ can be viewed as a bilinear functional on $M_*\times N(H)$.

From the other side, $M\bar\otimes B(K)\subset B(H\otimes K)$. We will use explicit coordinate form of this algebra.
For $x,y\in K$, denote by $\omega_{xy}$ the functional $\omega_{xy}(b)=\langle bx,y\rangle$, $b\in B(K)$.
For $a\in M$, $b\in B(K)$ we have with any $x,y\in H$, $x',y'\in K$:
\begin{align*}
(a\otimes b)(\mu_{xy},\omega_{x'y'}) &= \mu_{xy}(a)\omega_{x'y'}(b)
= \langle a x,y\rangle\langle b x',y'\rangle
\\& = \langle ax\otimes bx',y\otimes y'\rangle  = \langle (a\otimes b)(x\otimes x'),y\otimes y'\rangle.
\end{align*}
By continuity it follows that any $U\in M\bar\otimes B(K)\subset B(H\otimes K)$ acts as
\beq\label{U-on-HK}
U( \mu_{xy},\omega_{x'y'}) = \langle U (x\otimes x'),y\otimes y'\rangle.
\eeq

\begin{definition}
Let $\pi$ be a representation of $M_*$ on a Hilbert space $K$. An operator $U\in M\bar\otimes B(K)$ is called a {\it generator} of $\pi$ if
\beq\label{U-pi}
U(\mu,\omega) = \omega(\pi(\mu))
\eeq
for every $\mu\in M_*$, $\omega\in N(K)$.
\end{definition}

The original definition of a generator of a representation, where $U$ is supposed to be a partial isometry, is given in \cite[\S\,1.5]{enock}.

In fact, as shows the proposition below, the representations with generator are nothing else but completely bounded representations, as defined in the operator space theory \cite{effros}. In \cite{kirchberg}, the term ``bounded'' was used instead. When we speak of $M_*$ as an operator space, we always suppose the predual structure on it, that is the quotient structure in the corresponding $B(H)_*$ factored by $M_\perp$ (see \cite[4.2.2]{effros}).

\bprop\label{cb=U}
Let $\pi$ be a representation of $M_*$ on a Hilbert space $K$. Then $\pi$ is completely bounded if and only if there exists $U\in M\bar\otimes B(K)=(M_*\hat\otimes_{\rm op}N(K))^*$ such that $U(\mu,\omega) = \omega(\pi(\mu))$ for every $\mu\in M_*$, $\omega\in N(K)$. In this case $\|U\|=\|\pi\|_{cb}$.
\eprop
\bpr
There is \cite[7.1.5]{effros} a natural isometric isomorphism $\lambda$ of the spaces $CB(M_*,B(K))$ and
$M\bar\otimes B(K)$,
given by $\omega\big(\lambda(U)(\mu)\big) = U(\mu,\omega)$ for $U\in M\bar\otimes B(K)$ and $\mu\in M_*$, $\omega\in N(K)$. If $\pi\in CB(M_*,B(K))$, there is $U\in M\bar\otimes B(K)$ such that $\pi=\lambda(U)$, and vice versa. By definition of $\lambda$, this means exactly that $U(\mu,\omega) = \omega(\pi(\mu))$, and in this case $\|\pi\|_{cb}=\|U\|$.
\epr

If $\pi$ is a representation of $M_*$ on a Hilbert space $K$ with a basis $(f_\a)$, denote by $\pi_{\a\beta}\in M$ the linear functional on $M_*$ defined by $\pi_{\a\beta}(\mu)=\langle \pi(\mu)f_\beta,f_\a\rangle$, $\mu\in M_*$; in other words, $\pi_{\a\beta}(\mu) = \omega_{f_\beta,f_\a}(\pi(\mu))$.

We give next a temporary definition of a standard representation (this term is justified by Theorem \ref{continuous-is-standard} below). In the Theorem \ref{standard=generator} it is proved that a standard representation is in fact a non-degenerate representation with a unitary generator. There is another term for the same object: a {\it unitary representation}, and below we will switch to this latter term.

\begin{definition}\label{definition-standard}
Call a representation $\pi$ of $M_*$ on a Hilbert space $K$ {\it standard} if in some basis of $K$,
\beq\label{def-standard}
\sum_\gamma \pi_{\gamma\a}^*\cdot\pi_{\gamma\beta} = \sum_\gamma \pi_{\a\gamma}\cdot\pi_{\beta\gamma}^*
 = \begin{cases} 1, &\a=\beta\\ 0, &\a\ne\beta\end{cases}
\eeq
for every $\a,\beta$, the series converging absolutely in the $M_*$-weak topology of $M$.
\end{definition}

From the Theorem \ref{standard=generator} below it will follow that this definition does not depend, in fact, on the choice of a basis. This implies also that a standard representation is automatically non-degenerate.

\blm\label{mu-series}
Let $M$ act on a Hilbert space $H$ with a basis $(e_\a)$. Then for every $x,y\in H$ and every $a,b\in M$,
$$
(ab)(\mu_{xy}) = \sum_{\a} \mu_{e_\a y}(a) \mu_{xe_\a}(b),
$$
the series converging absolutely.
\elm
\bpr This is an immediate calculation:
\begin{align*}
(ab)(\mu_{xy}) &= \langle abx,y\rangle = \langle bx, a^*y\rangle
 = \sum_\a \langle bx,e_\a\rangle \langle e_\a,a^*y\rangle
\\& = \sum_\a \langle bx,e_\a\rangle \langle ae_\a,y\rangle
 = \sum_{\a} \mu_{e_\a y}(a) \mu_{xe_\a}(b).
\end{align*}
\epr

\blm\label{U-scalar-product}
Let $M\subset B(H)$, and let $\pi: M_*\to B(K)$ be a representation of $M_*$. Let $U\in M\bar\otimes B(K)$ be such that $U(\mu,\omega)=\omega(\pi(\mu))$ for every $\mu\in M_*$, $\omega\in N(K)$. Let $(e_\a), (f_\xi)$ be bases in $H$ and $K$ respectively. Then for any $x,y\in H$
$$
\langle U(x\otimes f_\a),U(y\otimes f_\beta)\rangle  = \sum_{\eta} ( \pi_{\eta\beta}^* \pi_{\eta\a})( \mu_{x,y}),
$$
$$
\langle U^*(x\otimes f_\a),U^*(y\otimes f_\beta)\rangle  = \sum_{\eta} (\pi_{\beta\eta}\pi_{\a\eta}^*) ( \mu_{x,y}),
$$
the series converging absolutely.
\elm
\bpr
With \eqref{U-on-HK}, $U$ satisfies for every $x,y\in H$ and every $\a,\beta$:
\begin{align}\label{def-U}
\langle U(x\otimes f_\a),(y\otimes f_\beta)\rangle &= U( \mu_{x,y},\omega_{f_\a,f_\beta})
 =\omega_{f_\a,f_\beta}( \pi(\mu_{x,y}) )
 = \pi_{\beta\a}( \mu_{x,y}),\\
\langle U^*(x\otimes f_\a),(y\otimes f_\beta)\rangle &= \overline{\langle U(y\otimes f_\beta),x\otimes f_\a\rangle}
 = \overline{\pi_{\a\beta}( \mu_{y,x})}.                                                              \notag
\end{align}
Decomposing the scalar product, we get the following absolutely converging series:
\begin{align}
\langle U(x\otimes f_\a),U(y\otimes f_\beta)\rangle&=\sum_{\xi,\eta} \langle U(x\otimes f_\a),e_\xi\otimes f_\eta\rangle
 \langle e_\xi\otimes f_\eta, U(y\otimes f_\beta)\rangle                                             \notag
\\&= \sum_{\xi,\eta} \pi_{\eta\a}( \mu_{x,e_\xi}) \overline{\pi_{\eta\beta}( \mu_{y,e_\xi})}.
\end{align}

After a simple transform:
$$
\overline{\pi_{\eta\beta}( \mu_{y,e_\xi})} = \overline{\langle \pi_{\eta\beta}\,y,e_\xi\rangle}
= \langle \pi_{\eta\beta}^*e_\xi,y\rangle= \pi_{\eta\beta}^*(\mu_{e_\xi,y}),
$$
we get:
$$
\langle U(x\otimes f_\a),U(y\otimes f_\beta)\rangle
 = \sum_{\xi,\eta} \pi_{\eta\a}( \mu_{x,e_\xi}) \pi_{\eta\beta}^*(\mu_{e_\xi,y}).
$$
With the lemma \ref{mu-series} this gives
$$
\langle U(x\otimes f_\a),U(y\otimes f_\beta)\rangle
 = \sum_{\eta} ( \pi_{\eta\beta}^* \pi_{\eta\a})( \mu_{x,y}).
$$

For $U^*$, we have similarly:
\begin{align*}
\langle U^*(x\otimes f_\a),U^*(y\otimes f_\beta)\rangle &
 = \sum_{\xi,\eta} \overline{\pi_{\a\eta}( \mu_{e_\xi,x})} \pi_{\beta\eta}( \mu_{e_\xi,y})
= \sum_{\xi,\eta} \pi_{\a\eta}^*( \mu_{x,e_\xi}) \pi_{\beta\eta}( \mu_{e_\xi,y})
\\&= \sum_{\eta} \pi_{\beta\eta}\pi_{\a\eta}^* ( \mu_{x,y}).
\end{align*}
\epr

\btm\label{standard=generator}
A non-degenerate representation of $M_*$ has a unitary generator if and only if it is standard.
\etm
\bpr
Let $M\subset B(H)$, and let $\pi:M_*\to B(K)$ be a representation. Choose bases $(e_\a)$, $(f_\beta)$ in the spaces $H$, $K$ respectively.

Suppose that a generator $U$ exists and is unitary. We have then for every $x,y\in H$ and every $\a,\beta$, by lemma \ref{U-scalar-product}:
$$
\langle x,y\rangle  \delta_{\a\beta} = \langle x\otimes f_\a,y\otimes f_\beta\rangle
 = \langle U(x\otimes f_\a),U(y\otimes f_\beta)\rangle
 = \sum_{\eta} ( \pi_{\eta\beta}^* \pi_{\eta\a})( \mu_{x,y}).
$$
Since $\langle x,y\rangle  = \mu_{x,y} (1)$, we get the first equality in \eqref{def-standard} for $\mu=\mu_{x,y}$. The decomposition into an absolutely converging series $\mu=\sum \mu_{x_n,y_n}$ implies that it is valid also for every $\mu\in M_*$.

The adjoint operator $U^*$ is unitary as well, and we get:
\begin{align*}
\langle x,y\rangle  \delta_{\a\beta}
=\langle U^*(x\otimes f_\a),U^*(y\otimes f_\beta)\rangle
= \sum_{\eta} \pi_{\beta\eta}\pi_{\a\eta}^* ( \mu_{x,y})
\end{align*}
what implies the second equality in \eqref{def-standard}, so that $\pi$ is standard.

Conversely, let $\pi$ be standard. Then we can take \eqref{def-U} as a definition of $U(x\otimes f_\a)$, and reversing the calculations above, we see that
$$
\langle U(x\otimes f_\a),U(y\otimes f_\beta)\rangle =  \langle x,y\rangle  \delta_{\a\beta}.
$$
Extending $U$ by linearity to finite linear combinations $x=\sum x_\a\otimes f_\beta$, we have:
$$
\|Ux\|^2 = \langle Ux,Ux\rangle
 = \sum_{\a,\beta,\gamma,\zeta} \langle U(x_\a\otimes f_\beta),U(x_\gamma\otimes f_\zeta)\rangle
 = \sum_{\a,\beta,\gamma,\zeta} \langle x_\a\otimes f_\beta,x_\gamma\otimes f_\zeta\rangle
= \|x\|^2.
$$
Thus, $U$ is isometric and then extends to an isometry on $H\otimes K$. Further, $U$ has an adjoint operator $U^*$ satisfying \eqref{def-U}; again, $\pi$ being standard implies that $U^*$ is isometric, so $U$ is onto and as a consequence unitary.

Now we must show that $U$ satisfies \eqref{U-pi}. By definition, we have \eqref{U-pi} for $\mu=\mu_{xy}$ and $\omega=\omega_{f_\a,f_\beta}$. Since both parts in
$$
\langle U(x\otimes x'),(y\otimes y')\rangle
 =\omega_{x',y'}( \pi(\mu_{x,y}) )
$$
depend linearly and jointly continuously on $x',y'\in K$, we have this equality for all $x',y'$, i.e.\ \eqref{U-pi} holds for $\mu=\mu_{xy}$ and $\omega=\omega_{x'y'}$. Finally, the decomposition $\mu=\sum \mu_{x_n,y_n}$, $\omega=\sum \omega_{x'_n,y'_n}$ into absolutely converging series implies \eqref{U-pi} for all $\mu\in M_*$, $\omega \in N(H)$.

Now it remains to show that $U\in M\bar\otimes B(H)$, not just $U\in B(H\otimes K)$. This follows from the bicommutant theorem; we need to show that $U$ commutes with $(M\bar\otimes B(H))'= M'\bar\otimes \C$ \cite[11.2.16]{kadison}. Let $\phi\in M'$. Then for every $x,y\in H$ and every $a\in M$ we have $\mu_{\phi x,y}(a)=\langle a\phi x,y\rangle = \langle \phi ax,y\rangle = \langle ax,\phi^*y\rangle=\mu_{x,\phi^*y}$, so that $\mu_{\phi x,y} = \mu_{x,\phi^*y}$. Then
\begin{align*}
\langle U(\phi\otimes \id)(x\otimes x'), y\otimes y'\rangle
&= \langle U(\phi x\otimes x'), y\otimes y'\rangle
 = \omega_{x',y'}(\pi(\mu_{\phi x,y}))
 = \omega_{x',y'}(\pi(\mu_{x,\phi^*y}))
 \\&= \langle U(x\otimes x'), \phi^*y\otimes y'\rangle
 = \langle (\phi\otimes \id)U(x\otimes x'), y\otimes y'\rangle,
\end{align*}
so that $U$ commutes with $\phi\otimes \id$, and the theorem is proved.
\epr

A finite-dimensional version of this theorem was known long ago \cite[Proposition 1.5.7]{enock}.

From now on, we will use a known term {\it unitary representation} instead of {\it standard}.
We see that the property of being unitary does not depend on the particular choice of a basis.

\bcor\label{standard-is-cb}
Every unitary representation $\pi$ of $M_*$ is completely contractive.
\ecor
\bpr
By Theorem \ref{standard=generator}, $\pi$ has a unitary generator; by Proposition \ref{cb=U} $\|\pi\|_{cb}=\|U\|=1$.
\epr

From \cite[3.1.4]{enock} it follows:
\bcor\label{Kac-every-rep-unitary}
If $M$ is a Kac algebra, every non-degenerate representation of $M_*$ is unitary.
\ecor

\section{Representations of the measure algebra}\label{reps-M(G)}

In this section we prove (Theorem \ref{standard=continuous}) that if an irreducible representation of $M(G)$ is unitary, then it comes from a continuous representation of $G$, and vice versa. Thus, in the case of the measure algebra, the class of unitary representations is exactly the class of representations commonly used in harmonic analysis.

The central theorem \ref{standard=continuous} in the abelian case reduces to a known theorem of M.~Walter \cite{walter}: a character (thus, an irreducible representation) of $B(G)$ is unitary, in the definition above, if and only if it is the evaluation at a point of $G$. Recall that if $G$ is abelian, then $B(G)\simeq M(\hat G)$. In general, Theorem \ref{standard=continuous} is a dual analogue of the Walter's theorem.

In principle, the results of this section are consequences of \cite{kirchberg}. But we find it highly instructive to present here new proofs, based on the representation theory. This gives a clear intuitive understanding of the whole picture.

\bprop\label{rep-on-L1}
Let $\pi$ be an irreducible representation of $M(G)$. Then either $\pi|_{L_1(G)}\equiv0$, or $\pi$ is generated by a continuous unitary representation $\tilde\pi$ of $G$ by the classical integral formula:
\beq\label{pi-integral}
\pi(\mu)=\int_G \tilde\pi(t)d\mu(t),
\eeq
for every $\mu\in M(G)$.
\eprop
\bpr
$L_1(G)$ is an ideal in $M(G)$, so the closed subspace $H_1$ generated by $\pi(f)H$, $f\in L_1(G)$, is invariant under $\pi$. Thus either $H_1=\{0\}$ or $H_1=H$. In the first case $\pi|_{L_1(G)}\equiv0$. In the second case we have a non-degenerate representation of $L_1(G)$, and one shows, as usual \cite[13.3]{dixmier}, that there is a continuous representation $\tilde\pi$ of $G$ such that \eqref{pi-integral} holds for all $\mu\in L_1(G)$.

It is known that a representation of a *-algebra is uniquely extended from a *-ideal on which it is non-degenerate \cite[11.1.12]{palmer}, so \eqref{pi-integral} holds for all $\mu\in M(G)$ as well.
\epr

\btm\label{continuous-is-standard}
Let $\bar\pi$ be a continuous unitary representation of \,$G$ on a Hilbert space $H$, and let $\pi$ be the corresponding representation of $M(G)$. Then $\pi$ is unitary.
\etm
\bpr
Let $(e_\a)_{\a\in \cal A}$ be a basis of $H$ and let $\pi_{\a\beta}(\mu) = \langle \pi(\mu) e_\beta,e_\a\rangle$. If $\mu=\delta_t$ is the probability measure in $G$ concentrated at a point $t$, we will also write $\pi_{\a\beta}(t):=\pi_{\a\beta}(\delta_t)$. Then $\pi_{\a\beta}$ may be considered as a continuous function on $G$.

To prove that $\pi$ is unitary, we need to show that for every $\mu\in M(G)$,
\begin{align}\label{mu-delta-alpha-beta}
\delta_{\a\beta} \int_G 1d\mu = \sum_\gamma \int_G \pi_{\gamma\a}(t) \overline{\pi_{\gamma\beta}(t)} d\mu
\end{align}
(since for continuous functions, the usual multiplication and conjugation coincide with those in $M(G)^*$). This corresponds to the first series in \eqref{def-standard}, interchanging $\a$ and $\beta$. The second series is reduced to the first one in the following way. One checks that $\pi_{\a\gamma}(t)=\overline{\pi_{\gamma\a}(t^{-1})} = \pi^*_{\gamma\a}(t^{-1})$. For $\mu\in M(G)$, let $\check\mu\in M(G)$ be defined by $\int f(t^{-1})d\mu(t)=\int f(t)d\check\mu(t)$. Then we have $\pi_{\a\gamma}(\mu)=\pi_{\gamma\a}^*(\check\mu)$. Once \eqref{mu-delta-alpha-beta} is proved, we will have then
$$
\sum_\gamma \int_G \pi_{\a\gamma} \pi_{\beta\gamma}^* d\mu =
\sum_\gamma \int_G \pi_{\gamma\a}^* \pi_{\gamma\beta} d\check\mu = \delta_{\a\beta}\int 1d\check\mu = \delta_{\a\beta} d\mu.
$$

The equality \eqref{mu-delta-alpha-beta} is sufficient to prove for a positive measure $\mu$.
For every $t\in G$, using the Kronecker symbol $\delta_{\a\beta}$, we have:
\begin{align}\label{delta-on-group}
\delta_{\a\beta}&=\langle e_\a,e_\beta\rangle = \langle \bar\pi(t)e_\a,\bar\pi(t)e_\beta\rangle = \sum_\gamma\langle \bar\pi(t)e_\a,e_\gamma\rangle \langle e_\gamma,\bar\pi(t)e_\beta\rangle
\\& = \sum_\gamma \pi_{\gamma\a}(t) \overline{\pi_{\gamma\beta}(t)}
 = \sum_{\gamma} (\pi_{\gamma\beta}^* \pi_{\gamma\a}) (\delta_t).                   \notag
\end{align}
Thus, pointwise (on $\delta_t$) we always have the equality \eqref{mu-delta-alpha-beta}. To prove that \eqref{mu-delta-alpha-beta} holds for any Radon measure (which we may assume to be positive), we need to use in further reasoning the continuity of $\pi$.

Denote
\beq\label{f-alpha-beta}
f_{\a\beta}(t) = \sum_\gamma |\pi_{\gamma\a}(t) \overline{\pi_{\gamma\beta}(t)}|.
\eeq
This series converges, since
\begin{align*}
f_{\a\beta}(t)
& = \sum_\gamma |\langle \bar\pi(t)e_\a,e_\gamma\rangle|\,|\langle e_\gamma,\bar\pi(t)e_\beta\rangle|
 \\&\le \Big(\sum_\gamma |\langle \bar\pi(t)e_\a,e_\gamma\rangle|^2\Big)^{1/2}
 \,\Big(\sum_\zeta |\langle e_\zeta,\bar\pi(t)e_\beta\rangle|^2\Big)^{1/2}
\\&=\|\bar\pi(t)e_\a\|\,\|\bar\pi(t)e_\beta\| = 1.
\end{align*}
If $\a=\beta$, then $f_{\a\a}\equiv1$. Let us show that $f_{\a\beta}$ is a continuous function, whatever are $\a$ and $\beta$. For $s,t\in G$,
\begin{align*}
|f_{\a\beta}(t)-f_{\a\beta}(s)|
&= \sum_\gamma \Big( |\pi_{\gamma\a}(t) \pi_{\gamma\beta}(t)| - |\pi_{\gamma\a}(s) \pi_{\gamma\beta}(s)|\Big)
\\&= \sum_\gamma \Big( |\pi_{\gamma\a}(t) \pi_{\gamma\beta}(t)| - |\pi_{\gamma\a}(t) \pi_{\gamma\beta}(s)| + |\pi_{\gamma\a}(t) \pi_{\gamma\beta}(s)|- |\pi_{\gamma\a}(s) \pi_{\gamma\beta}(s)|\Big)
\\&\le \sum_\gamma \Big( |\pi_{\gamma\a}(t)|\cdot| \pi_{\gamma\beta}(t) -\pi_{\gamma\beta}(s)| + |\pi_{\gamma\beta}(s)|\cdot |\pi_{\gamma\a}(t)- \pi_{\gamma\a}(s)|\Big)
\\&\le \Big(\sum_\zeta |\pi_{\zeta\a}(t)|^2 \sum_\gamma | \pi_{\gamma\beta}(t)\! -\!\pi_{\gamma\beta}(s)|^2\Big)^{1/2} + \Big( \sum_\zeta |\pi_{\zeta\beta}(s)|^2 \sum_\gamma |\pi_{\gamma\a}(t)\!- \!\pi_{\gamma\a}(s)|^2\Big)^{1/2}
\\&= \Big(\sum_\zeta |\langle \bar\pi(t)e_\a,e_\zeta\rangle|^2 \sum_\gamma | \langle \bar\pi(t)e_\beta,e_\gamma\rangle -\langle \bar\pi(s)e_\beta,e_\gamma\rangle |^2\Big)^{1/2}
\\&+ \Big( \sum_\zeta |\langle \bar\pi(s)e_\beta,e_\zeta\rangle|^2 \cdot \sum_\gamma |\langle \bar\pi(t)e_\a,e_\gamma\rangle - \langle \bar\pi(s)e_\a,e_\gamma\rangle |^2\Big)^{1/2}
\\&= \|\bar\pi(t)e_\a\| \| \bar\pi(t)e_\beta - \bar\pi(s)e_\beta\| + \|\bar\pi(s)e_\beta\| \cdot \| \bar\pi(t)e_\a - \bar\pi(s)e_\a\|
\\&= \| \bar\pi(t)e_\beta - \bar\pi(s)e_\beta\| + \| \bar\pi(t)e_\a - \bar\pi(s)e_\a\|.
\end{align*}
As $s\to t$, this tends to zero since $\bar\pi$ is continuous in the strong-operator topology (as every continuous unitary representation).
Thus, the series \eqref{f-alpha-beta} of positive continuous functions converges to a continuous function; by Dini's theorem, it converges uniformly on every compact subset of $G$.

We can assume that $\mu(G)=1$. For every $\e>0$, there is a compact set $F\subset G$ such that $\mu(G\setminus F)<\e$ \cite[14.1]{HR}. Now on $F$, the series converges absolutely and uniformly to a constant function $\delta_{\a\beta}$ (see \eqref{delta-on-group}), so we immediately get \eqref{mu-delta-alpha-beta}.
\epr

\bprop\label{L1-invariant}
The annihilator of $L_1(G)$ is a two-sided ideal in $M(G)^*$.
\eprop
\bpr
Denote $M=M(G)^*= C_0(G)^{**}$. This space can be viewed as the enveloping von Neumann algebra of $C_0(G)$. By the corresponding universality property, every representation $\rho$ of $C_0(G)$ is extended to a normal representation $\bar\rho$ of $M$, so that $\bar\rho(M)$ is the von Neumann algebra generated by $\rho(C_0(G))$. Let $\rho$ be the action on $L_2(G)$ by pointwise multiplication. Then $\bar\rho(M)=L_\infty(G)$.

From the other side, consider $L_1(G)$ as a closed subspace of $M(G)$, and let $Z=L_1(G)^\perp$ be the annihilator of $L_1(G)$ in $M$. Then $L_1(G)^*=M/Z$. This space is also isomorphic to $L_\infty(G)$, so we get the quotient map $q: M\to L^\infty(G)$.  Then $\bar\rho$, $q$ are both weakly continuous and identical on $C_0(G)$; since $C_0(G)$ is weakly dense in $M$, it follows that $\bar\rho=q$. Thus, $\ker \bar\rho=\ker q=L_1(G)^\perp$, so this is a two-sided ideal in $M$.
\epr

\btm\label{standard=continuous}
An irreducible representation $\pi$ of $M(G)$ is unitary if and only if it is generated by a continuous representation $\tilde\pi$ of $G$ by the integral formula \eqref{pi-integral}.
\etm
\bpr
With Theorem \ref{continuous-is-standard}, we need to prove one implication only. Let the restriction of $\pi$ to $G$ be discontinuous. Then, by Proposition \ref{rep-on-L1}, $\pi|_{L_1(G)}\equiv0$. By Proposition \ref{L1-invariant}, then $(\pi_{\gamma\a}^* \pi_{\gamma\beta})|_{L_1(G)}\equiv0$ for every $\a,\beta,\gamma$; thus
$$
\sum_{\gamma} (\pi_{\gamma\a}^* \pi_{\gamma\beta})(f) = 0
$$
for any $f\in L_1(G)$. This is different from $f(1)=\int f$ if $\int f\ne0$, so $\pi$ is not unitary.
\epr

\bprop
Every representation of $M(G)$ is completely bounded.
\eprop
\bpr
We always consider $M(G)$ with the operator space structure as the dual of $C_0(G)$, or, what is the same, as the predual of $C_0(G)^{**}$. Since $C_0(G)$ is a commutative \cst-algebra, its natural operator structure is the minimal one. Then on $M(G)$ we have the maximal operator space structure. And in this case, every bounded linear operator on $M(G)$ is completely bounded.
\epr

\section{The absolutely continuous ideal}

\begin{definition}\label{def-Mhat}
Representations of $M_*$ which are not unitary we will call non-unitary. Let $M_*^\times\subset M_*$ be the common kernel of all {\it irreducible} non-unitary representations. If there are none, let $M_*^\times=M_*$. This is a two-sided *-ideal in $M_*$, which is called the {\it absolutely continuous ideal} of $M_*$. With the structure inherited from $M_*$, $M_*^\times$ is a Banach *-algebra and an operator space.
\end{definition}

\brem\label{nondegen=contained-in-closure}
Every non-degenerate representation of $M_*^\times$ extends uniquely to $M_*$ \cite[11.1.12]{palmer}. It is easy to show that $M_*$ is mapped into the weak closure of the image of $M_*^\times$, i.e. to the von Neumann algebra generated by $M_*^\times$. Conversely, if $\phi_1,\phi_2: M_*\to N$ are two *-homomorphisms to a von Neumann algebra $N$ which agree on $M_*^\times$ and are such that $\phi_i(M_*)$ is contained in the weak closure of $\phi_i(M_*^\times)$ for $i=1,2$, then $\phi_1=\phi_2$.
\erem

In the case when $M_*=M(G)$, the algebra $M_*^\times$ has been studied by J.~Taylor under the notation $L^{1/2}(G)$ \cite{taylor}. He has proved that $L^{1/2}(G) \ne L_1(G)$ if $G$ is non-discrete. However, $L^{1/2}(G)$ and $L_1(G)$ have the same $*$-representations, what motivates our term: this ideal is a means of `recovering' the absolutely continuous measures without knowing the Haar measure.

The main property of $M_*^\times$ is that every representation which is non-degenerate on it must be unitary. For the duality construction we need more: that unitary representations do not vanish on it. To guarantee this, we exclude all degenerate cases by the following definition:

\begin{definition}
Let $I^0$ be the weakly closed ideal in $M$ generated by $(M_*^\times)^\perp$, that is by the annihilator of $M_*^\times$. Set $M_*^0=M_*^\times$ if $I^0\ne M$ and $M_*^0=\{0\}$ otherwise.
\end{definition}

\bprop\label{unitary-implies-non-degenerate}
If $M_*^0\ne\{0\}$ then every unitary representation of $M_*$ is non-degenerate on $M_*^0$.
\eprop
\bpr
Let $\pi: M\to B(H)$ be a unitary representation of $M_*$. First, it is nonzero on $M_*^0$: otherwise we would have $\pi_{\a\beta}\in I^0$ for all its coefficients, and by \eqref{def-standard} this would imply $1\in I^0$, what is not true by assumption.

Next suppose that $\pi$ is degenerate on $M_*^0$. Let $L\subset H$ be the null subspace of $\pi(M_*^0)$. Then $L$ is $\pi$-invariant, and $\pi|_L$ is also unitary \cite[1.5.4 (ii)]{enock} and vanishes on $M_*^0$, what is impossible.
\epr

The main property of the ideal $M_*^0$ is presented in the next theorem. Recall that every (completely bounded) irreducible representation which does not vanish on $M_*^0$ is unitary by definition.

The main property of the ideal $M_*^0$ is presented in the next theorem.

\btm\label{nondegenerate-is-standard}
Every completely bounded representation of $M_*$ which is nondegenerate on $M_*^0$ is unitary.
\etm
\bpr
By Proposition \ref{cb=U}, there exists $U\in M\bar\otimes B(H)$ such that $U(\mu,\omega) = \omega(\pi(\mu))$, and all we need is to prove that $U$ is unitary.

Let $M$ be realized on a Hilbert space $K$. We will need several times the following representation. Fix $x\in K$, $\xi\in H$ and bases $(e_\a)\subset K$, $(f_\beta)\subset H$. Then
$$
\langle U(x\otimes \xi),e_\a\otimes f_\beta\rangle = \langle \pi(\mu_{x e_\a}) \,\xi,f_\beta\rangle,
$$
so that
\begin{align}\label{U-series}
U(x\otimes \xi)
 &= \sum_{\a,\beta} \langle U(x\otimes \xi),e_\a\otimes f_\beta\rangle\, e_\a\otimes f_\beta
 = \sum_{\a,\beta} \langle \pi(\mu_{x e_\a}) \,\xi,f_\beta\rangle \,e_\a\otimes f_\beta                      \notag
\\ &= \sum_\a e_\a\otimes \Big( \sum_\beta \langle \pi(\mu_{x e_\a})\, \xi,f_\beta\rangle \;f_\beta \Big)
= \sum_\a e_\a\otimes \pi(\mu_{x e_\a})\, \xi
\end{align}
(convergence is everywhere in the Hilbert space norm).

For subspaces $E,F\subset K$, let $M_{E,F}\subset M_*$ denote the closed subalgebra generated by $\mu_{xy}$ with $x\in E$, $y\in F$. Denote also $M^0_{E,F} = M_{E,F}\cap M^0_*$. These subalgebras are not supposed to be self-adjoint. By $\langle M_{E,F}\rangle^*$, $\langle M_{E,F}^0\rangle^*$ we denote the closed $*$-subalgebras generated by $M_{E,F}$ and $M_{E,F}^0$ respectively.

We can suppose that $M$ is realized in its standard form. Then \cite[1.2.8]{enock} the involution on $M_*$ is given by $\mu_{xy}^* = \mu_{Jx,Jy}$, with an antilinear bijective isometry $J:K\to K$. For $J$-invariant subspaces $E$ and $F$, we have then $\langle M_{E,F}\rangle^*=M_{E,F}$ and $\langle M_{E,F}^0\rangle^*=M_{E,F}^0$.

\blm\label{invar-subspace-U}
A closed subspace $E\otimes L\subset K\otimes H$ is $U$-invariant if and only if $\pi(M_{E,K})L\subset L$ and $\pi(M_{E,E^\perp})L=\{0\}$.
\elm
\bpr
Choose a base $(e_\a)_{\a\in A}$ in $K$ such that the subset $(e_\a)_{\a\in A_1}$ is a base for $E$.
If $E\otimes L$ is $U$-invariant, then, by \eqref{U-series}, for every $x\in E$, $\xi\in L$ by orthogonality $\pi(\mu_{x e_\a})\, \xi\in L$ for all $\a$, and moreover $\pi(\mu_{x e_\a})\, \xi=0$ for $\a\notin A_1$. It follows that $\pi(M_{E,K})L\subset L$ and $\pi(M_{E,E^\perp})L=\{0\}$.

Conversely, if $\pi(M_{E,E^\perp})$ vanishes on $L$, then for all $x\in E$, $\xi\in L$ the sum in \eqref{U-series} reduces to $\a\in A_1$ only, and if $L$ is invariant under $\pi(M_{E,K})$, then moreover $U(x\otimes \xi)\in E\otimes L$.
\epr

\blm\label{invar-subspace-F-UF}
For every separable subspace $F\subset K\otimes H$ there exist closed separable subspaces $E\subset K$, $L\subset H$ such that $F\cup UF \subset E\otimes L$ and $(J\otimes {\rm id})F\subset E\otimes L$.
\elm
\bpr
Since $F$ is separable, so is $V=\overline{F+UF}$. Pick a sequence $(x_n)$ dense in $V$ and orthonormal bases $(e_\a)\subset K$, $(f_\beta)\subset H$. Every $x_n$ is contained in $\overline{\rm lin}\{ e_\a\otimes f_\beta: \a\in A_n, \beta\in B_n\}$ with countable $A_n$, $B_n$. Then for $E_0=\overline{\rm lin}\{ e_\a: \a\in \cup A_n\}$ and $L = \overline{\rm lin}\{ f_\beta: \beta\in \cup B_n\}$ we have $F\cup UF \subset E_0\otimes L$. The statement now holds with $E=\overline{E_0+JE_0}$.
\epr

\blm
Every $v\in K\otimes H$ can be embedded into a $U$-invariant separable subspace $E\otimes L$ such that $JE\subset E$ and $L$ is essential for $\pi(M_{E,K}^0)|_L$.
\elm
\bpr
Construct separable subspaces $E_k$, $L_k$ by induction as follows. Let $E_1\otimes L_1$ be any separable subspace containing $v$. Suppose now that $E_{k-1}$, $L_{k-1}$ are constructed for some $k\ge2$. Since $H$ is essential for $M_*^0$, there are sequences $\mu^{(k)}_n\in M_*^0$, $\xi^{(k)}_n\in H$ such that ${\rm lin}\{ \pi(\mu^{(k)}_n)\xi^{(k)}_n\}$ is dense in $L_{k-1}$. Since $M$ is in the standard form, every $\mu^{(k)}_n$ can be represented as $\mu^{(k)}_n=\mu_{x^{(k)}_n,y^{(k)}_n}$ with $x^{(k)}_n,y^{(k)}_n\in K$. Set $E'_k=\overline{ E_{k-1}+{\rm lin}\{ x^{(k)}_n: n\in\N \}}$, $L'_k=\overline{L_{k-1}+ {\rm lin}\{ \xi^{(k)}_n: n\in\N \}}$. Then $\mu^{(k)}_n \in M_{E'_k,K}^0$ and $\xi^{(k)}_n\in L'_k$ for all $n$. By Lemma \ref{invar-subspace-F-UF} there are separable subspaces $E_k$, $L_k$ such that $E'_k\otimes L'_k\cup U(E'_k\otimes L'_k) \subset E_k\otimes L_k$ and $JE'_k\subset E_k$.

Set  $E=\overline{\cup E_k}$ and $L=\overline{\cup L_k}$, then $E\otimes L = \overline{ \cup (E_k\otimes L_k)}$ since $E_k$, $L_k$ are increasing. We have $U(E_k\otimes L_k)\subset E_{k+1}\otimes L_{k+1}$, $JE_k\subset E_{k+1}$ for all $k$, what implies the invariance required in the statement. Moreover, by construction the set $\{\pi(\mu)\xi: \mu\in M^0_{E,K},\xi\in L\}$ is dense in $L$.
\epr

{\it Proof of the theorem.} Take now any $v\in K\otimes H$. Let $E\otimes L\subset K\otimes H$ be $U$-invariant and separable, such that $v\in E\otimes L$, $JE\subset E$ and $L$ is essential for $M_{E,K}^0$. Together with $E$, the subspace $E^\perp$ is also $J$-invariant (since $J$ is isometric), what implies that $M_{E,E}$, $M_{E,E}^0$, $M_{E,E^\perp}$, $M_{E,K}$ and $M_{E,K}^0$ are self-adjoint.

Let $r_L:B(H)\to B(L)$ be the reduction onto $L$, then we have a *-representation $\rho=r_L\circ\pi: M_{E,K}\to B(L)$ (such that $\rho(M_{E,E^\perp})=0$). Let $\cal A$ be the \cst-algebra generated by $\rho(M_{E,K}^0)$ in $B(L)$.
Since $M_{E,E}$ is separable, so are $\rho(M_{E,E})=\rho(M_{E,K})\supset \rho(M_{E,K}^0)$ and $\cal A$.

The identity representation of $\cal A$ is then decomposed into a direct integral of irreducible representations \cite[8.5.2]{dixmier}: there exist a set $P$ equipped with a probability measure $\beta$; an integrable field of Hilbert spaces $\Gamma\subset \{(H_p)_{p\in P}\}$; a field of representations $\sigma_p: {\cal A}\to B(H_p)$, $p\in P$, where every $\sigma_p$ is irreducible; and an isometric isomorphism $V:L\to \Gamma=\int_P^\oplus H_p d\beta(p)$ such that: if $\xi\in L$ and $V(\xi)=\int_P \xi_p \,d\beta(p)$,
then for every $a\in \cal A$ we have $V(a\xi)=\int_P \sigma_p(a)\xi_p \,d\beta(p)$.

For every $\mu\in M_{E,K}^0$ and $\xi\in L$, denoting $V\xi=\int\xi_p$, we have
\beq\label{rho-equal-int}
V(\rho(\mu)\xi)=\int_P \rho_p(\mu)\xi_p \,d\beta(p)
\eeq
with irreducible representations $\rho_p=\sigma_p\circ\rho$ of $M^0_{E,K}$.

Moreover, we can extend $\rho_p$ uniquely and irreducibly to $M_{E,K}$, still denoting this extension by $\rho_p$.
We have $\rho(M_{E,K})\subset {\cal A}''$ \cite[14.1.10]{kadison}, so every $V\circ\rho(\mu)$ with $\mu\in M_{E,K}$ is decomposable, and one verifies that the formula \eqref{rho-equal-int} holds in fact for all $\mu\in M_{E,K}$.

Let $[\cdot]$ denote the norm closure in $B(H)$. For every $p$, $\sigma_p\circ r_L$ is lifted from $[\pi(M^0_{E,K})]$ to an irreducible representation $\tilde\sigma_p$ of $[\pi(M_*^0)]$, probably on a bigger space $\tilde H_p\supset H_p$ \cite[2.10]{dixmier}. Next, $\tilde\sigma_p\circ\pi$ is extended uniquely from $M_*^0$ to an irreducible representation $\tilde \rho_p$ of $M_*$. On $M^0_{E,K}$, $\tilde\rho_p|_{H_p}=\rho_p$; since $\rho_p$ is irreducible (and cyclic) on the ideal $M_{E,K}^0\subset M_{E,K}$, its extension to $M_{E,K}$ is unique, and we have $\tilde\rho_p|_{H_p}=\rho_p$ on $M_{E,K}$.

In particular, $H_p$ is invariant under $\tilde\rho_p(M_{E,K})$, and $\tilde\rho_p(M_{E,E^\perp})|_{H_p}=\rho_p(M_{E,E^\perp})=\{0\}$ (this follows from $\rho(M_{E,E^\perp})=\{0\}$ and \eqref{rho-equal-int}).
Being irreducible and nonzero on $M_*^0$, $\tilde\rho_p$ is unitary, with a unitary generator $U_p\in M\otimes B(\tilde H_p)$. By the reasoning above, we can apply Lemma \ref{invar-subspace-U} and conclude that $E\otimes H_p$ is $U_p$-invariant.

$\int(E\otimes H_p)_{p\in P}$ is also a field of Hilbert spaces, isomorphic to $E\otimes L$ under the isomorphism $\tilde V=\id\otimes V$. Let  $(e_\a)_{\a\in A}$ be a base in $K$ such that the (countable) subset $(e_\a)_{\a\in A_1}$ is a base for $E$. We have then for $x\in E$, $\xi\in L$ that $\mu_{x e_\a}\in M_{E,K}$ and so $\pi(\mu_{x e_\a})\xi=\rho(\mu_{x e_\a})\xi$. Denote $V\xi=\int \xi_p$; we have
\begin{align}
\tilde V(U(x\otimes \xi)) &= \tilde V\Big(\sum_{\a\in A_1} e_\a\otimes \pi(\mu_{x e_\a})\, \xi\Big)         \notag
\\&= \sum_{\a\in A_1} e_\a\otimes V(\rho(\mu_{x e_\a})\, \xi)
 = \sum_{\a\in A_1} e_\a\otimes \int \big(\rho_p(\mu_{x e_\a})\, \xi_p\big)                                 \notag
\\&= \sum_{\a\in A_1} \int e_\a\otimes \big(\rho_p(\mu_{x e_\a})\, \xi_p\big).                              \label{series-ea-int}
\end{align}
The last series converges in the Hilbert norm of $\int(E\otimes H_p)$.

From the other hand, for every $p$ we have a formula similar to \eqref{U-series}: if $x\in E$, $\xi_p\in H_p$, then
$$
U_p(x\otimes \xi_p) = \sum_{\a\in A_1} e_\a\otimes \rho_p(\mu_{x e_\a})\, \xi_p,
$$
so the series in \eqref{series-ea-int} converges pointwise to $\int U_p(x\otimes \xi_p)$. Both imply convergence in measure in the following sense: denote  $\phi_{\a p}=e_\a\otimes\rho_p(\mu_{x e_\a})\, \xi_p$, then for every $\e>0$
\begin{align*}
&\beta\{ p: \|\big(\tilde V(U(x\otimes \xi))\big)_p - \sum_{\a\in B} \int \phi_{\a p}\|\ge\e\} \to0,\\
&\beta\{ p: \|U_p(x\otimes \xi_p) - \sum_{\a\in B} \int \phi_{\a p}\|\ge\e\} \to0
\end{align*}
as finite set of indices $B\subset A_1$ increases (the reasoning for real-valued functions applies verbatim). It follows that $\big(\tilde V(U(x\otimes \xi))\big)_p = U_p(x\otimes \xi_p)$ almost everywhere, that is
$$
\tilde V(U(x\otimes \xi)) = \int U_p(x\otimes \xi_p) = \Big(\int U_p\Big)(\tilde V(x\otimes \xi)).
$$
It follows that $(U_p)$, or strictly speaking $(U_p|_{E\otimes H_p})$, is a measurable field of operators on $\int E\otimes H_p$, and $\tilde V U = \int U_p\tilde V$. Since $\int U_p$ is unitary, so is $U$ on $E\otimes L$.

As the initial  vector $v$ was arbitrary, we get that $U$ is unitary on $H\otimes K$, what proves the theorem.
\epr

\bcor
If $M_*^0$ is nonzero, then a completely bounded representation of $M_*$ is unitary if and only if it is non-degenerate on $M_*^0$.
\ecor

We give next several examples: first some degenerate ones, and then we consider
the case of a locally compact group. The dual algebra will be defined via $C^*(M_*^0)$, this is why we mention this \cst-algebra in every case.

\bex\label{ex-dual-is-C}
On $M=L^\infty(\R)$ with the usual structure, introduce a new coinvolution $\kappa$ as $\kappa(\phi)=\phi$. Then $M$ is again a coinvolutive Hopf--von Neumann algebra. On its predual, which is $L_1(\R)$, we have the usual convolution product. The involu\-{tion ${}^\circ$} on $L_1(\R)$ is the pointwise conjugation: for $f\in L_1(\R)$, $\phi\in L_\infty(\R)$ we have
$$
\int \phi f^\circ = \overline{\int (\kappa\phi)^* f } = \int \phi(t) \overline{f(t)}dt,
$$
so $f^\circ(t)=\overline{f(t)}$.

Every $s\in \R$ defines by the usual formula a character $\chi_s$ of $L_1(\R)$, and every character is of this form. A character is involutive if
$$
\chi_s(f^\circ) = \int e^{ist} \overline{f(t)}dt = \overline{\chi_s(f)} = \int e^{-ist}\overline{f(t)}dt,
$$
i.e. if $s=0$. For $s=0$ the trivial character $\chi_0$ is of course unitary. Thus, there are no non-unitary characters, so $M_*^0=M_*^\times=M_*$. And it is now obvious that $C^*(M_*)=\C$. This shows in particular that the canonical map of $M_*^0$ to $C^*(M_*^0)$ might not be injective.
\eex

\bex
Modifying the previous example a little bit, one can show that an algebra ``arising'' from one locally compact group can be driven to the group algebra of another group by our functor.

Namely, let $M=L^\infty(\R^2)$. Change the coinvolution to be $\kappa(\phi)(s,t)=\phi(s,-t)$. Then on $M_*=L_1(\R^2)$ we get as before the usual convolution and the involution $f^\circ(s,t)=\overline{f(s,-t)}$. The {\it involutive\/} characters are described as $\chi_{\a}(f) = \int e^{i\a t} f(s,t)dsdt$, $\a\in\R$ and are all unitary, thus we have again $M_*^0=M_*^\times=M_*$ but $C^*(M_*^0)=C_0(\R)$.
\eex

\bex\label{ex-dual-is-zero}
It may happen that $M_*^\times=\{0\}$. This means that non-unitary irreducible representations separate points of $M_*$.

Let $S$ be a topological semigroup with identity. Call a map $V:S\to B(H)$ a representation of $S$ if $V(xy)=V(x)V(y)$ and $\|V(x)\|\le1$ for all $x,y\in S$. If $x$ is invertible, it follows that $V(x)$ must be unitary. A representation will be called continuous if it is a continuous map from $S$ to $B(H)$ in the weak operator topology. One can show \cite[p.~55]{kirchberg}
that there is a von Neumann algebra $W^*(S)$ such that its normal representations correspond bijectively to continuous representations of $S$. Moreover, there is a continuous map with a dense image $F:S\to W^*(S)$ (where $W^*(S)$ is considered with the weak topology), and $W^*(S)$ has a canonical structure of a Hopf--von Neumann algebra. If $S$ is a locally compact group, then $W^*(S)$ is exactly the Ernest group algebra. One can introduce also coinvolution on $W^*(S)$ so that it induces the pointwise involution on the predual $B(S)=(W^*(S))_*$  (this is the only detail not mentioned in \cite{kirchberg}). With this structure, $W^*(S)$ is a coinvolutive Hopf--von Neumann algebra.

By definition, $F(x)$ is an involutive character of $B(S)$ for every $x\in S$. If $F(x)$ is a unitary element in $W^*(S)$, then $V(x)$ must be unitary for every representation $V$. Now it is easy to provide an example of $S$ such that $F(x)$ is not unitary for all $x$ except the identity. For example, if $S=[0,+\infty)$ (with addition), then for every $\a<0$ we have a non-unitary character $V(x)=e^{\a x}$, $x\in S$, so we have $V\in B(S)$. If $x\ne0$, then $|V(x)|=|e^{\a x}|\ne1$, so $F(x)$ is not unitary for any $x\ne0$, so this is a non-unitary character of $B(S)$. For $x=0$, the character $F(0)$ is of course unitary.
From the other side, $F(S\setminus \{0\})$ separates points of $B(S)$ by continuity, so $M_*^\times=M_*^0=\{0\}$.

We see also on this example that the set of representations of $C^*(M_*^0)$ may be strictly less than the set of unitary representations of $M_*$.
\eex

\bprop Let $A$ be a Banach *-algebra, and let $I\subset A$ be a
two-sided *-ideal. Let $\phi: A\to C^*(A)$ be the canonical map.
Then $C^*(I)$ is isomorphic to the closure of $\phi(I)$ in
$C^*(A)$. \eprop \bpr Let $\gamma_A$, $\gamma_I$ be the maximal
\cst-seminorms on $A$ and on $I$ respectively. Clearly
$\gamma_A|_I\le \gamma_I$. From the other side, let $\pi:I\to
B(H)$ be a representation of $I$ such that
$\|\pi(x)\|=\gamma_I(x)$ for all $x\in I$. Then
\cite[11.1.12]{palmer} $\pi$ can be extended to a representation
$\bar\pi:A\to B(H)$, and we will have
$\|\pi(x)\|\le\|\bar\pi(x)\|\le\gamma_A(x)$ for all $x\in I$. It
follows that $\gamma_I=\gamma_A|_I$.

Since $C^*(I)$ is the completion of $I$ with respect to $\gamma_I$
and $C^*(A)$ is complete, from the isometry $\phi(I)\simeq
(I,\gamma_I)$ it follows the isometry in question.
\epr

\bprop Let $A$ be a Banach *-algebra, and let $I\subset A$ be a
two-sided *-ideal. Let $\phi: A\to W^*(A)$ be the canonical map.
Then $W^*(I)$ is isomorphic to the weak closure of $\phi(I)$ in
$W^*(A)$. \eprop \bpr Denote by $[\phi(I)]$ this weak closure; it
is a von Neumann algebra. Let $\pi$ be a representation of $I$. It
has a unique extension $\bar\pi$ to $A$, and $\bar\pi$ has a
unique normal lifting to $W^*(A)$. Its restriction to $[\phi(I)]$
is a normal lifting of $\pi$, in the same sense that it is normal
and equals to $\pi$ when composed with $\phi$.

To show that this lifting is unique, suppose that $\rho$, $\sigma$ are two distinct normal representations of $[\phi(I)]$ both lifting $\pi$. We can assume that $\pi$ is non-degenerate, then so are $\rho$ and $\sigma$. By \cite[11.1.12]{palmer}, there are unique extensions $\bar\rho$, $\bar\sigma$ to $W^*(A)$ (since $[\phi(I)]$ is an ideal in $W^*(A)$). They are obviously given by $\bar\rho(x)=\rho(px)$ and $\bar\sigma(x)=\sigma(px)$, where $p\in W^*(A)$ is a central projection such that $pW^*(A)=[\phi(I)]$. Thus, these extensions are normal. Since $\phi(A)$ is weakly dense in $W^*(A)$, $\bar\rho$ and $\bar\sigma$ have different restrictions to $\phi(A)$. But it would mean that $\pi$ has two different extensions to $A$, what is impossible.

Thus, $[\phi(I)]$ has the universal property of the von Neumann
envelope of $I$, so it is isomorphic to $W^*(I)$. \epr

Let $A$ be a Banach *-algebra. For a set $X\subset A$, let $h(X)$
be the set of irreducible representations of $A$ which vanish on
$X$.

\bprop\label{envelope-equal} Let $A$ be a Banach *-algebra. Let
$B$, $B'$ be ideals in $A$ such that $h(B)=h(B')$. Then
$C^*(B)\simeq C^*(B')$. \eprop \bpr Let $\phi:A\to C^*(A)$ be the
canonical map, and for $X\subset A$ let $[X]$ denote the closure
of $\phi(X)$ in $C^*(A)$. As we have proved above,
$C^*(B)\simeq[B]$ and $C^*(B')\simeq[B']$. It is easy to see that
$[B], [B']$ are ideals in $C^*(A)$.

Let $\Omega$ denote the space of irreducible representations of
$C^*(A)$. For $Y\subset \Omega$, set $k(Y)=\cap\{\ker\pi: \pi\in
Y\}$. In $C^*(A)$, as in every \cst-algebra, $I=k(h(I))$ for every
closed ideal $I$ \cite[2.9.7]{dixmier}. From assumptions it
follows that $h([B])=h([B'])$.

Thus, $C^*(B)\simeq[B]= k(h([B]))=k(h([B']))=[B']\simeq C^*(B')$.
\epr

\bprop\label{classical-envelopes}
If $M_*=M(G)$ or $M_*=L_1(G)$, then $C^*(M_*^0)=C^*(L_1(G))=C^*(G)$. If $M_*=B(G)$ or $M_*=A(G)$, then $C^*(M_*^0)=C_0(G)$.
\eprop
\bpr
For $M_*=M(G)$, let $I=M(G)^\times$ be the absolutely continuous  ideal. It contains obviously $L_1(G)$ but is strictly
larger than $L_1(G)$ if $G$ is non-discrete, see \cite{taylor}. Since $L_1(G)^\perp$ is a proper ideal in $M(G)^*$, $M_*^0=I$.
By the results of Sect\-{ion \ref{reps-M(G)}}, $h(I)=h(L_1(G))$. Then, by Proposition \ref{envelope-equal}, $C^*(I)=C^*(L_1(G))=C^*(G)$. If $M_*=L_1(G)$, then $M_*^0=M_*$ and also
$C^*(M_*^0)=C^*(L_1(G))=C^*(G)$.

In the second case, the algebra $M_*=B(G)$ is commutative, so its
irreducible representations are characters. A character is unitary
if and only if it is unitary in $M=W^*(G)$. It was proved by
M.~Walter \cite{walter} that a character of $B(G)$ is unitary if
and only if it is the evaluation at a point $t\in G$, and if and
only if it does not vanish on $A(G)$. Thus, $B(G)^\times\supset A(G)$ and $h(B(G)^\times) = h(A(G))$. Since $A(G)^\perp$ is a proper ideal in $W^*(G)=B(G)^*$, it follows that $B(G)^0=B(G)^\times$
and we can again apply Proposition \ref{envelope-equal}, to
conclude that $C^*(B(G)^0)=C^*(A(G))$. By Proposition
\ref{envelope-A(G)}, $C^*(A(G))=C_0(G)$. And finally, if
$M_*=A(G)$, then by Corollary \ref{Kac-every-rep-unitary}
$M_*=M_*^0$, and $C^*(A(G)^0) = C^*(A(G))=C_0(G)$.
\epr

This proposition already establishes a duality in the group case:
\beq\label{square-cstar}
\xymatrix @R=1.pc @C=2.pc
 {
 C_0(G)\ar@{|->}[rrr]^{*}
 &
 & &
 M(G) \ar@{|->}[dd]^{C^*(M_*^0)}
 \\
 & & &
 \\
 B(G) \ar@{|->}[uu]^{C^*(M_*^0)}
 & & &
 C^*(G)\ar@{|->}[lll]_{*}
 }
\eeq In the next section, we develop this construction to a
general framework of Hopf-von Neumann algebras, and in Section
\ref{dual-cstar}, to the category of coinvolutive \cst-bialgebras.

\section{The dual Hopf-von Neumann algebra}\label{dual-von-neumann}

Throughout the section $M$ will denote a coinvolutive Hopf--von
Neumann algebra.

\let\tilde\widetilde
\begin{definition}
Set $\Mhat=W^*(M_*^0)$ be the enveloping von Neumann algebra of $M_*^0$. There is a canonical map $\Phi: M_*^0\to \Mhat$, maybe not injective. It has a unique extension to $M_*$, still having range in $\Mhat$ (see discussion after Definition \ref{def-Mhat}), which we denote also $\Phi: M_*\to \Mhat$. By the universality property, every representation of $M_*^0$ is lifted to a unique normal representation of $\Mhat$.
\end{definition}

By definition, $\Mhat$ is a von Neumann algebra. Below we define on it a structure of a coinvolutive Hopf--von Neumann algebra.

The coinvolution on $\Mhat$ is given by the composition with $\kappa$ (cf. \cite[1.6]{enock}):

\bprop\label{dual-kappa}
Let $\tilde\kappa:M_*\to M_*$ be defined by $\tilde\kappa(\mu)=\mu\circ\kappa$ for $\mu\in M_*$. Then $\tilde\kappa$ is a *-antihomomorphism, maps $M_*^0$ to itself, and is uniquely lifted from $M_*^0$ to a normal *-antihomomorphism $\hat\kappa$ of $\Mhat$, which is a coinvolution on $\Mhat$.
\eprop
\bpr
It is proved, for example, in \cite[1.6]{enock} that $\tilde\kappa$ is a *-antimohomorphism and $\tilde\kappa^2=\id$.

To prove that $M_*^0$ is invariant under $\tilde\kappa$, we should show that $\pi(\tilde\kappa(M_*^0))=0$ for every non-unitary irreducible representation $\pi$ of $M_*$. Fix $\pi$ and set $\tau_{\a\beta}=\kappa\pi_{\beta\a}$ for every $\a,\beta$. Then, by \cite[1.4.2]{enock}, these are coefficients of a representation of $M_*$, since
$$
\kappa(\tau_{\a\beta}^*) = \kappa(\kappa\pi_{\beta\a}^*)=\kappa(\pi_{\a\beta})=\tau_{\beta\a},
$$
\begin{align*}
\Delta(\tau_{\a\beta}) &= \Delta(\kappa\pi_{\beta\a})=\theta(\kappa\otimes\kappa)\Delta(\pi_{\beta\a})
 = \theta(\kappa\otimes\kappa)\Big(\sum_\gamma \pi_{\beta\gamma}\otimes\pi_{\gamma\a}\Big)
 \\&= \sum_\gamma \kappa(\pi_{\gamma\a})\otimes \kappa(\pi_{\beta\gamma})
 = \sum_\gamma \tau_{\a\gamma}\otimes \tau_{\gamma\beta},
\end{align*}
the series converging ultraweakly.
If $\pi$ acts on a Hilbert space $H$ and $J:H\to \bar H$ is the anti-isomorphism onto the conjugate space $\bar H$, then we can view $\tau$ as acting on $\bar H$ with $\tau(\mu)=J\pi(\tilde\kappa\mu^*)J^{-1}$. Since $\pi$ is by assumption irreducible and $\tilde\kappa(M_*)=M_*$, this implies that $\tau$ is irreducible.

If $\tau$ were unitary, then by \eqref{def-standard} we would have, since $\kappa$ is ultraweakly continuous:
$$
\sum_\gamma \pi_{\gamma\a}^*\pi_{\gamma\beta} = \sum_\gamma \kappa(\tau_{\a\gamma})^*\kappa(\tau_{\beta\gamma})
 = \sum_\gamma \kappa(\tau_{\beta\gamma}\tau_{\a\gamma}^*) =  \kappa\Big( \sum_\gamma \tau_{\beta\gamma}\tau_{\a\gamma}^*\Big)
  = \kappa(\delta_{\a\beta})=\delta_{\a\beta},
$$
and similarly $\sum \pi_{\a\gamma}\pi_{\beta\gamma}^*=\delta_{\a\beta}$. Thus, $\pi$ would be also a unitary representation, what is not true by assumption. This implies that $\tau$ is non-unitary.

Now $\tau_{\a\beta}(\mu)=0$ for all $\mu\in M_*^0$, so $\pi_{\beta\a}(\tilde\kappa\mu)=0$ and $\pi(\tilde\kappa(M_*^0))=0$. Since $\pi$ was arbitrary, this implies that $M_*^0$ is invariant under $\tilde\kappa$.

By universality, $\tilde\kappa$ is extended to a *-antihomomorphism $\hat\kappa$ of $\Mhat$. Similarly to \cite[1.6.6]{enock}, one proves  that $\theta(\hat\kappa\otimes\hat\kappa)\hat\Delta=\hat\Delta\tilde\kappa$.
\epr

\bprop\label{Phi-is-cb} The canonical map $\Phi: M_*\to \Mhat$ is
completely contractive. \eprop \bpr
Let $\pi$ be the direct sum of all irreducible representations of $M_*$ which are nonzero on $M_*^0$. They are in bijection with the irreducible representations of $M_*^0$. Let us consider $\pi$ as a map to the \cst-algebra $A$ generated by $\pi(M_*)$.
By the theory above and \cite{enock}, $\pi$ has a unitary
generator, so it is completely contractive.

Now $C^*(M_*^0)$ is isometrically isomorphic to the closure of
$\pi(M_*^0)$ in $A$ (denote this closure by $I$) and to the norm
closure of $\Phi(M_*^0)$ in $\Mhat=W^*(M_*^0)$. Let $\rho:I\to
\Mhat$ be the latter isomorphism. Clearly $I$ is an ideal in $A$,
so $\rho$ is extended to a unique representation $\bar\rho$ of $A$
having range in $\Mhat$. As every homomorphism of \cst-algebras,
$\bar\rho$ is completely contractive. Now $\Phi
=\bar\rho\circ\pi$, so $\Phi$ is also completely contractive. \epr

\bcor\label{Phi-is-unitary-if-separable}
If $\hat M$ is nonzero, then $\Phi$ is a unitary representation of $M_*$.
\ecor
\bpr
This follows from Proposition \ref{Phi-is-cb} and Theorem \ref{nondegenerate-is-standard}, since $\Phi$ is by assumption non-degenerate on $M_*^0$.
\epr

\brem\label{PhixPhi=DeltaPhi}
Recall that on the set of representations of $M_*$, there is a so
called Kronecker product \cite[1.4.3]{enock}, denoted by $\times$.
We have therefore a representation $\Phi\times\Phi$ of $M_*$ whose
coefficients are by definition $(\Phi\times\Phi)_*(\omega\otimes
\upsilon) = \Phi_*(\omega)\cdot\Phi_*(\upsilon)$, for
$\omega,\upsilon\in \Mhat_*$. Denote by $\hat\Delta: \Mhat\to \Mhat\bar\otimes\Mhat$ the normal lifting of $\Phi\times\Phi$ to
$\Mhat$, so that $\hat\Delta\Phi=\Phi\times\Phi$ on $M_*^0$.
In fact, $\Phi\times\Phi$ and $\hat\Delta\Phi$ are (by \cite[1.5.5]{enock} and Remark \ref{unitary-implies-non-degenerate} respectively) both non-degenerate (or null) on $M_*^0$ and equal on it, so the equality $\hat\Delta\Phi=\Phi\times\Phi$ holds on $M_*$ and not just on $M^0_*$.
\erem

$\Phi$ has the adjoint map $\Phi^*: \hat M^*\to M$. Let $\hat\Phi:
\hat M_*\to M$ be the restriction of $\Phi^*$ onto $\hat M_*\equiv
(\hat M)_*$. Sometimes we will write the corresponding algebra $M$
as an index, so that we have $\Phi_M: M_*\to \hat M$ and
$\hat\Phi_M: \hat M_*\to M$. Note that $\hat\Phi_M$ is always
injective, because $\Phi_M(M_*)$ is weakly dense in $\hat M$.

\bprop\label{dual-comultiplication}
$\hat\Delta$ is a comultiplication on $\hat M$.
\eprop
\bpr
The statement is trivial if $\Mhat=\{0\}$, so we assume further that $\Mhat$ is nonzero.
First, we need to prove that $\hat\Delta$ is coassociative: $(\id\otimes\hat\Delta)\hat\Delta=(\hat\Delta\otimes\id)\hat\Delta$. We have $\hat\Delta\Phi=\Phi\times\Phi$ (see Remark \ref{PhixPhi=DeltaPhi}). This implies that $\hat\Phi\hat\Delta_* = (\Phi\times\Phi)_*: (\hat M\bar\otimes \hat M)_*\to M$. Take now $\mu\in M_*$, $\omega_j\in \hat M_*$, $j=1,2,3$. First,
\let\o\omega
$$
((\Phi\times\Phi)_*(\o_1\otimes\o_2))(\mu) = (\o_1\otimes\o_2)((\Phi\times\Phi)(\mu))=(\hat\Phi(\o_1)\cdot\hat\Phi(\o_2))(\mu).
$$
Using this identity, we can calculate the following:
\begin{align*}
(\o_1\otimes\o_2\otimes\o_3)(\id\otimes\hat\Delta)\hat\Delta(\Phi(\mu))
&= (\o_1\otimes(\hat\Delta_*(\o_2\otimes\o_3)))(\Phi\times\Phi)(\mu)
\\&= (\hat\Phi(\o_1)\cdot \hat\Phi(\hat\Delta_*(\o_2\otimes\o_3)))(\mu)
\\&= (\hat\Phi(\o_1)\cdot \hat\Phi(\o_2)\cdot\hat\Phi(\o_3))(\mu)
\end{align*}
It is easy to see that we come to the same result, starting with $(\hat\Delta\otimes\id)\hat\Delta$ instead, so associativity holds.

Next, $\hat\Delta$ is injective. To show this, take first $\omega\in\hat M_*$ and $\mu\in M_*$.  Note that the unit of $M$ is a *-homomorphism from $M_*$ to $\C$, so it has a lifting to $\hat M$. We can thus consider $1$ as an element of $\hat M_*$. Now
$$
(\o\otimes1)(\hat\Delta(\Phi\mu)) = (\o\otimes1)(\Phi\times\Phi(\mu))=(\hat\Phi(\o)\cdot\hat\Phi(1))(\mu)=\o(\Phi(\mu)).
$$
It follows that $(\o\otimes1)(\hat\Delta(x))=\o(x)$ for all $x\in\hat M$. This shows that $\hat\Delta$ is injective.

By Corollary \ref{Phi-is-unitary-if-separable}, $\Phi$ is unitary, and by \cite[1.5.5]{enock}, $\Phi\times\Phi=\hat\Delta\Phi$ is unitary too. It is then non-degenerate, what implies $\hat\Delta(1)=1\otimes 1$ (cf. \cite[1.6.5]{enock}), and we conclude that $\hat\Delta$ is a comultiplication on $\Mhat$.
\epr

\begin{definition}
The algebra $\Mhat$ with the comultiplication and coinvolution introduced above is called the {\it dual coinvolutive Hopf-von Neumann algebra} of $M$.
\end{definition}

The structure of $\Mhat$ gives rise to a structure of a Banach *-algebra on $\Mhat_*$. From the equality $\hat\Delta\Phi=\Phi\times\Phi$ it follows immediately that $\hat\Phi$ is a *-homomorphism.

\bcor\label{Phi-is-standard}
If $\hat M$ is nonzero, then $\hat\Phi$ is a unitary representation of\/ $\hat M_*$.
\ecor
\bpr
Let $U\in M\bar\otimes\hat M$ be the unitary generator of $\Phi$. Then $\theta U\in \hat M\bar\otimes M$ is the unitary generator for $\hat\Phi$, where $\theta$ is the flip: $\theta U(\omega,\mu)=U(\mu,\omega)=\Phi(\mu)(\omega)=\mu(\hat\Phi(\omega))$ for all $\mu\in M_*$, $\omega \in\hat M_*$.
\epr

\bprop\label{dual-morphism}
Let $\phi: M\to N$ be a morphism of coinvolutive Hopf--von Neumann algebras. Then there is a dual morphism $\widehat\phi:\hat N\to \hat M$ such that $\hat\phi\circ\Phi_N=\Phi_M\circ\phi_*$ on $N_*^0$. If $\hat N\ne\{0\}$, then the equality holds on $N_*$.
\eprop
\bpr
The statement is trivial if $N_*^0=\{0\}$ so we can assume that $N_*^0\ne\{0\}$.
By definition, $\phi$ is ultraweakly continuous, so it has a pre-adjoint $\phi_*:N_*\to M_*$. Since $\phi$ is a coalgebra morphism, $\phi_*$ is a *-homomorphism. Consider $\psi=\Phi_M\circ\phi_*: N_*^0\to \Mhat$. This is a *-homomorphism, so it is lifted to a normal *-homomorphism $\hat\phi: \hat N\to \hat M$, such that $\hat\phi\,\Phi_N=\psi$ (on $N_*^0$).

If $\Mhat=\{0\}$ then trivially $\hat\phi\,\Phi_N=\psi$ on $N_*$. Suppose further that $\Mhat\ne\{0\}$. Then $\hat\Phi_M$ is unitary, and (let $[\cdot]$ denote the weak closure) $[\hat\Phi_M(\Mhat_*)]=[\hat\Phi_M(\Mhat^0_*)]$. By weak continuity, then $\phi(\hat\Phi_M(\Mhat_*))\subset [\phi(\hat\Phi_M(\Mhat^0_*))]$. This means that $\phi\circ\hat\Phi_M$ is non-degenerate on $\Mhat^0_*$, and by Theorem \ref{nondegenerate-is-standard} it is unitary. Similarly to Corollary \ref{Phi-is-standard}, then $\Phi_M\phi_*=(\phi\hat\Phi_M)_*$ is unitary too, and as a consequence non-degenerate on $N_*^0$. From the other side, since $\Phi_N$ is unitary and $\hat\phi$ weakly continuous, $\hat\phi\,\Phi_N(N_*)\subset \hat\phi\,\Phi_N(N_*^0)$, so that $\hat\phi\,\Phi_N$ is also non-degenerate on $N^0_*$. Two representations of $N_*$ which are equal and non-degenerate on $N^0_*$ must be equal everywhere, so in fact the equality holds on $N_*$.

It remains to prove that $\hat\phi$ is a coalgebra morphism: $\Delta_{\Mhat} \, \widehat\phi = (\widehat\phi\otimes\widehat\phi)\Delta_{\widehat N}$. Since $\Delta_{\Mhat}$ is ultraweakly continuous, this equality is enough to check on $\Phi_N(N_*^0)$. Moreover, to check an equality in $\Mhat\bar\otimes\Mhat$ where $\hat\Delta$ takes its values, it is enough to consider evaluations on $x\otimes y$, with $x,y\in \Mhat_*$. We have, with any $\nu\in N_*^0$:
\begin{align}\label{nu-calculations}
\Delta_{\Mhat} \, \widehat\phi (\Phi_N(\nu))(x\otimes y) &= \widehat\phi(\Phi_N(\nu))(xy) = \Phi_M(\phi_*(\nu))(xy)
=\phi_*(\nu)(\hat\Phi_M(x)\hat\Phi_M(y))
\notag\\&=
\nu\Big(\phi(\hat\Phi_M (x))\,\phi(\hat\Phi_M (y))\Big).
\end{align}
By definition, $\Delta_{\widehat N}\Phi_N = \Phi_N\times\Phi_N$, so from the other side:
\begin{align*}
(\widehat\phi\otimes\widehat\phi)\Delta_{\widehat N}(\Phi_N(\nu))(x\otimes y)
&= (\Phi_N\times\Phi_N)(\nu)(\hat\phi_* x\otimes \hat\phi_*y)
= \nu\big((\Phi_N)_*(\hat\phi_* x)\cdot (\Phi_N)_*(\hat\phi_*y)\big).
\end{align*}
Recalling that $\hat\phi\,\Phi_N = \Phi_M\phi_*$, we arrive at the required equality.

Also, $\kappa_{\widehat M}\widehat\phi = \widehat\phi\kappa_{\widehat N}$ should hold. Again, it is sufficient to consider $\xi=\Phi_N(\nu)$, $\nu \in N_*^0$. Then, with the definition of $\hat\kappa$ from Proposition \ref{dual-kappa}:
$$
\kappa_{\widehat M}\widehat\phi(\Phi_N(\nu)) = \kappa_{\widehat M}\Phi_M(\phi_*(\nu)) = \Phi_M(\tilde\kappa_M(\phi_*(\nu))).
$$
For every $x\in \hat M_*$ we have further
\begin{align*}
\Phi_M(\tilde\kappa_M(\phi_*(\nu)))(x) &= \tilde\kappa_M(\phi_*(\nu))(\hat\Phi_M(x))
 = \nu\big(\phi(\kappa_M\hat\Phi_M(x))\big) = \nu\big(\kappa_N(\phi(\hat\Phi_M(x)))\big)
\\&= \tilde\kappa_N(\nu)\big((\Phi_N)_*(\hat\phi_*(x))\big)
=\hat\phi(\Phi_N\tilde\kappa_N(\nu))(x)
= \hat\phi \hat\kappa_N(\Phi_N(\nu))(x),
\end{align*}
as required.
\epr

\brem\label{hat-is-functor}
If we have two morphisms $\phi:M\to N$ and $\psi: N\to L$, then on $L_*^0$ we have:
$$
\hat{\psi\circ\phi}\circ\Phi_L = \Phi_M\circ (\psi\circ\phi)_* = \Phi_M\circ\phi_*\circ\psi_*,
$$
$$
\hat\phi\circ\hat\psi \circ\Phi_L = \hat\phi \circ \Phi_N \circ \psi_*.
$$
If $\hat N\ne\{0\}$, then we can use the fact that $\hat\phi \circ \Phi_N=\Phi_M\circ\phi_*$ on $N_*$ and conclude that both displayed lines are equal, so that $\hat{\psi\circ\phi}=\hat\phi\circ\hat\psi$. If $\hat N=\{0\}$, then necessarily $\hat\phi\circ\hat\psi=0$ but it might happen that $\hat{\psi\circ\phi}\ne0$.
\erem

We have finally the following duality theorem:

\btm\label{von-neumann-duality}
Let $\cal H_0$ be the category of coinvolutive Hopf--von Neumann algebras, and let $\cal H$ be the full subcategory of $\cal H_0$ which has as its objects all $M\in \cal H_0$ such that $M\simeq \hat{\Mhat}$; such $M$ is called reflexive. Then $\hat{\;}:\cal H\to \cal H$ is a contravariant duality functor, i.e. $\hat{\Mhat}\simeq M$ for every $M\in \cal H$. Moreover, there is a faithful contravariant functor $\cal A:\cal{LCG}\to\cal H$, where $\cal{LCG}$ is the category of locally compact groups, such that $\hat{\cal A(G)} = \cal A(\hat G)$ for every abelian group $G\in \cal{LCG}$.
\etm
\bpr
The subcategory $\cal H$ is full by Proposition \ref{dual-morphism}, and $\hat{\;\;}$ is by definition a duality functor.

For every $G\in \cal{LCG}$, set $\cal A(G) = C_0(G)^{**}$. It is known that this is a faithful contravariant functor from $\cal {LCG}$ to $\cal H_0$ (or it follows from \cite[5.1.4]{enock} and Remark \ref{hat-is-functor}). And it follows from Proposition \ref{classical-envelopes} that $C_0(G)^{**}$ is reflexive, i.e. $\cal A(G)\in \cal H$.
\epr

\section{Particular cases}

In this section we consider first the Kac algebras and relate their duals in the sens above to the usual duals (Proposition \ref{dual-of-Kac-is-reflexive}). Next, Propositions \ref{D-exists} and \ref{DE=id} establish relations between $M$ and its second dual. This allows to describe, in Theorems \ref{commutative-is-group} and \ref{cocommutative-is-group}, commutative and co-commutative algebras in the range of our functor. These are nothing but the algebras $C_0(G)^{**}$ or $C^*(G)^{**}$, where $G$ is a locally compact group. In particular, this describes commutative and co-commutative reflexive algebras. In \cite[1.7]{kirchberg} there is a similar characterization, but only for reflexive algebras, and with some additional assumptions in the co-commutative case.

The following proposition has been already proved by E.~Kirchberg \cite{kirchberg}. In his definition, the ``unitary dual'' $W^*U(M_*)$ of $M$ is a coinvolutive Hopf--von Neumann algebra whose normal representations are in bijection with {\it all} unitary representations of $M$. If $M$ is a Kac algebra, then our $\hat M$ coincides with the unitary dual of Kirchberg. In general, they may differ. In the example \ref{ex-dual-is-zero} $\hat M=\{0\}$ while $W^*U(M_*)=\C$.

\bprop\label{dual-of-Kac-is-reflexive}
Let $N$ be a Kac algebra. Let $L$ be the Kac dual of $N$ and let $M=\hat N$. Then $M=W^*(N_*)$ and $\hat L=\hat M=\hat{\hat N}=W^*(L_*)$. In particular, $\hat{\hat M}\simeq M$.
\eprop
\bpr
Recall that a Kac algebra is supposed to be nonzero. It is known that every representation of $N_*$ has a unitary generator \cite[3.1.4]{enock}, so $N_*^0=N_*$. By definition, $M=\hat N$ is then $W^*(N_*)$, the von Neumann envelope of $N_*$.

Clearly $L$ is a quotient of $M$: $L=\phi(M)$, where $\phi$ is a morphism of coinvolutive Hopf--von Neumann algebras. Then $L_*=(\ker\phi)^\perp$ can be considered as an ideal of $M_*$. If we prove that $L_*=M_*^0$, it will follow that $\hat M=W^*(L_*)$.

From one side, if an irreducible representation $\pi$ is zero on $L_*$, then its coefficients are contained in $(L_*)^\perp=\ker\phi\subset M$. Since $\ker\phi$ is a proper ideal, $\pi$ is non-unitary.

It remains to show, from the other side, that every irreducible representation of $M_*$ is unitary if its restriction to $L_*$ is nonzero. In its turn, it is enough to show that the extension $\tau$ to $M_*$ of the universal representation $\Phi_L:L_*\to W:=W^*(L_*)$ is unitary.

It is known \cite[3.2.2(ii)]{enock} that $\Phi_L$ is quasi-equivalent to the regular representation $\Lambda:L_*\to N$ (based on the Haar weight), so that it can be expressed through $\Lambda$ as $(\Lambda\times \Phi_L)(\rho) = Z(\Lambda(\rho))$, $\rho\in L_*$, where $Z$ is an isomorphism of von Neumann algebras. One can express then the generator $U$ of $\Phi_L$ explicitly via the generator $W$ of $\Lambda$ (as in \cite[Theorem 3.1.4]{enock}): $U\otimes 1 = (1\otimes \theta)(W^*\otimes 1)(\theta\otimes 1)(1\otimes Z) W$ in $L\otimes W^*(L)\otimes N$, with the help of the flip $\theta$ on $W^*(L)\otimes L$.

On $N_*$, we have a representation $\kappa_L\circ\lambda:N_*\to L$ (which is unitary as every non-degenerate representation of $N_*$). It factors through $M=W^*(N_*)$: $\kappa_L\circ\lambda =\sigma\circ \Phi_N$ with some morphism $\sigma: M\to L$. Taking pre-adjoints in this equality, we get $(\kappa_L\circ\lambda)_* = (\Phi_N)_*\sigma_* = \hat\Phi_N\sigma_*$, with $\sigma_*:L_*\to M_*$ being an injection. From the other side, $(\kappa_L\circ\lambda)_* = \Lambda$ \cite[3.7.3]{enock}, and we see that $\hat\Phi_N$ is an extension of $\Lambda$ to $M_*$, both representations being unitary.

Now both $\hat\Phi_N\times \tau$ and $Z\hat\Phi_N$, as representations of $M_*$, extend the same non-degenerate representation $\Lambda\times \Phi_L = Z\Lambda$ of $L_*$. It follows that they are equal. Let $\tilde W$ be the generator of $\hat\Phi_N$. Then a formula similar to that of $U$ gives the generator $\tilde U$ of $\tau$: $\tilde U\otimes 1 := (1\otimes \theta)(\tilde W^*\otimes 1)(\theta\otimes 1)(1\otimes Z) \tilde W$, in $M\otimes W^*(L)\otimes N$, and this proves the proposition.
\epr

If $M$ is a finite-dimensional Kac algebra, then $W^*(M_*)=\hat M$ is also equal to the Kac dual of $M$ (and is also finite-dimensional), see \cite[6.6.9]{enock}. In particular, $\hat{\hat M}=M$.

If $M$ is infinite-dimensional, then its second dual is usually not equal to $M$. But, by Proposition \ref{dual-of-Kac-is-reflexive}, the first dual equals to the third dual, and so the first dual is a reflexive algebra in the sense of our duality.

The table below summarizes results on the duals of classical algebras, given by Propositions \ref{classical-envelopes} and \ref{dual-of-Kac-is-reflexive}.
$$
\begin{tabular}{|ll|}
\hline
\vrule height14pt width0pt\ Algebra\quad\null& Dual\\
\hline
\vrule height20pt width0pt$\left.\begin{matrix}L_\infty(G)\\C_0(G)^{**}\end{matrix}\right.$&$W^*(G)$\\
\hline
\vrule height26pt width0pt$\left.\begin{matrix}\cal L(G)\\W^*(G)
\end{matrix}\right.$&$C_0(G)^{**}$\\
\hline
\end{tabular}
$$
\vskip2pt

It is clear that not every coinvolutive Hopf--von Neumann algebra is reflexive. For instance, this is shown by the example \ref{ex-dual-is-zero} of a nontrivial algebra with $\widehat M=0$. Below we establish some relations between $M$ and $\hat{\hat M}$ in the general case.

\bprop\label{D-exists}
For every $M$, there is a morphism $D_M: \widehat{\Mhat}\to M$ such that $D_M\circ\Phi_{\Mhat}(x)=\hat\Phi_M(x)$ for all $x\in \Mhat_*^0$. If $\hat{\Mhat}\ne\{0\}$, then this equality holds actually on $\Mhat_*$.
\eprop
\bpr
Denote $N=\Mhat$. We know that there is a canonical map $\hat \Phi_M: N_*\to M$. In particular, $\hat\Phi_M|_{N_*^0}$ is a *-homomorphism, so it is extended to a normal homomorphism of its von Neumann envelope: $D_M:\widehat N\to M$. By definition, $D_M$ satisfies the equality in the statement for $x\in N_*^0=\Mhat_*^0$.

If $\Mhat_*^0=\{0\}$, then $D_M=0$ is by definition a morphism. Assume further that $\Mhat_*^0\ne\{0\}$. Then $\hat\Phi_M$ is unitary and non-degenerate on $\Mhat_*^0$. In the equality $D_M\circ\Phi_{\Mhat}=\hat\Phi_M$, we have then two representations, equal and non-degenerate on $\Mhat_*^0$; it follows that they are equal on $\Mhat_*$.
The preadjoint map $(\hat\Phi_M)_*: M_*\to \hat M$ is by definition $\Phi_M$, but from the equality above it is equal also to $(\Phi_N)_*\circ(D_M)_* = \hat\Phi_N\circ (D_M)_*$.

To prove that $D_M$ is a morphism of coinvolutive Hopf--von Neumann algebras, one should check the equality
\begin{align*}
\Delta_M \, D_M (x)(\mu\otimes \nu) = (D_M\otimes D_M)\Delta_{\widehat N}(x)(\mu\otimes \nu)
\end{align*}
for every $x\in \hat N$ and $\mu,\nu\in M_*$. By density, it is sufficient to consider $x=\Phi_N(y)$ with $y\in N_*^0$. Then we have:
\begin{align*}
\Delta_M \, D_M\Phi_N(y)(\mu\otimes \nu) &= \Delta_M \, \hat\Phi_M(y)(\mu\otimes \nu)=
\hat\Phi_M(y)(\mu\nu) = y(\Phi_M(\mu\nu))
\\&= y( \Phi_M(\mu)\Phi_M(\nu)) = y(\hat\Phi_N((D_M)_*(\mu))\cdot\hat\Phi_N((D_M)_*(\nu)))
\\&=\Phi_N(y)((D_M)_*(\mu)\cdot(D_M)_*(\nu)) = \Delta_{\widehat N}(x)((D_M)_*(\mu)\otimes(D_M)_*(\nu))
\\&= (D_M\otimes D_M)(\Delta_{\widehat N}(x))(\mu\otimes \nu).
\end{align*}

Also, $\kappa_M D_M= D_M\kappa_{\widehat N}$ should hold. Using the definition of the dual coinvolution in Proposition \ref{dual-kappa} (and its notations, so that $\tilde\kappa_M(\mu)=\mu\circ\kappa_M$), we have for every $y\in \hat M_*^0$, $\mu\in M_*$ (note that $\tilde\kappa_N(y)\in \Mhat_*^0$ by Proposition \ref{dual-kappa}):
\begin{align*}
\kappa_M D_M \Phi_N(y)(\mu) &= \kappa_M \hat\Phi_M(y)(\mu) = \hat\Phi_M(y)(\tilde\kappa_M(\mu))=
 y(\Phi_M(\tilde\kappa_M(\mu)))=
 \\&= y(\kappa_{\hat M}(\Phi_M(\mu)) = \tilde\kappa_N(y)(\Phi_M(\mu)) = \hat\Phi_M(\tilde\kappa_N(y))(\mu)=
 \\&=D_M \Phi_N\tilde\kappa_N(y)(\mu) =D_M \kappa_{\hat N} \Phi_N(y)(\mu).
\end{align*}
This shows that $\kappa_M D_M= D_M\kappa_{\widehat N}$ on $\Phi_N(N_*^0)$, and by density on the whole of $\hat{\hat M}$.
\epr

One should note that $D_M$ need not be neither injective nor surjective. As a first example, take the algebra $M$ from the example \ref{ex-dual-is-zero}. Then $\hat\Mhat=0$, so $D_M=0$. For a second example, take $M=L_\infty(G)$. Then $\hat\Mhat=C_0(G)^{**}$, and $D_M$ is a quotient map but is not injective.

However, if $M$ is a dual of another algebra, then $D_M$ is right invertible:

\bprop\label{DE=id}
If $M=\hat N$ for some $N$ and $\hat{\Mhat}\ne\{0\}$, then there is a morphism $E_M: M\to \widehat {\widehat M}$, such that $D_M\circ E_M=\id_M$. If $\hat{\hat {\hat M}}\ne\{0\}$ then $E_M$ is unital.
\eprop
\bpr
By the previous proposition, there is a morphism $D_N: \widehat{\hat N}\to N$ such that $D_N\circ\Phi_{\hat N}(x)=\hat\Phi_N(x)$ for all $x\in \hat N_*$. By Proposition \ref{dual-morphism}, there exists a dual morphism $E_M=\widehat D_N: \hat N=M\to (\widehat{\hat N})\widehat{\;} = \widehat{\widehat M}$.

By assumption $\Mhat=\hat{\hat N}\ne0$, and in this case it was proved in Proposition \ref{D-exists} that $\Phi_N = \hat\Phi_M\circ (D_N)_*$ on $N_*$. Moreover, since $\hat{\Mhat}\ne\{0\}$, we have $D_M\circ \Phi_{\hat M}=\hat\Phi_M$ on $M_*$.

By Proposition \ref{dual-morphism}, since $M\ne\{0\}$, we have $E_M\circ\Phi_N=\Phi_{\hat M}\circ (D_N)_*$ on $N_*$. Then
$$
D_M\circ(E_M\circ\Phi_N) = (D_M\circ \Phi_{\hat M})\circ (D_N)_* = \hat\Phi_M\circ (D_N)_* = \Phi_N.
$$
Thus, $D_M\circ E_M=\id_M$ on $\Phi_N(N_*)$. Since this latter is weakly dense in $M$, this equality holds everywhere.

This implies, in particular, that $D_M$ is surjective, and that $\hat\Phi_M(\hat M^0_*)$ is weakly dense in $M$. If $\hat{\hat {\hat M}}\ne\{0\}$ (in fact, even if $\hat{\hat M}\ne\{0\}$), then we can apply this reasoning to $\hat M$ and conclude that $\hat\Phi_{\hat M}(\hat {\Mhat}{}^{0}_*)$ is weakly dense in $\hat M$.

Suppose now that $\hat{\hat {\hat M}}\ne\{0\}$ but $E_M(1)=p\ne1$. Then $E_M(M)$ is contained in the weakly closed *-subalgebra $I:=p\hat{\hat M}p$.
Consider the representation $\pi=\Phi_M\circ E_{M*}:\hat{\hat M}_*\to \hat M$. Its coefficients for $\omega\in \hat M_*$ are: $\omega(\pi(\mu)) = \Phi_M\circ E_{M*}(\mu)(\omega) = \mu(E_M\circ\hat\Phi_M(\omega))$, so the space of coefficients is contained in the subalgebra $I$. Then the equality \eqref{def-standard} cannot hold, so $\pi$ is not unitary. From Propositions \ref{Phi-is-cb} and \ref{nondegenerate-is-standard} it follows that $\pi$ is degenerate on $\hat{\hat M}{}_*^0$, or equivalently $\pi(\hat{\hat M}{}_*)$ is not contained in the weak closure of $\pi(\hat{\hat M}{}_*^0)$. In particular, $\pi(\hat{\hat M}{}_*^0)$ is not weakly dense in $\Mhat$.

Consider now $\hat E_M:\hat{\hat {\hat M}}\to \hat M$. From one side, $\hat E_M\Phi_{\hat{\hat M}}=\Phi_ME_{M*}=\pi$ on $\hat{\hat M}{}_*^0$, and since $\hat E_M$ is weakly continuous, $\hat E_M(\hat{\hat {\hat M}})$ is contained in the closure of $\hat E_M(\Phi_{\hat{\hat M}}(\hat{\hat M}{}_*^0))= \pi(\hat{\hat M}{}_*^0)$, so $\hat E_M$ is not surjective. From the other side, $\hat E_M\hat D_M=(D_ME_M)\hat{\;}={\hat{\Id}}_M=\Id_{\Mhat}$ so $\hat E_M$ must be surjective. This contradiction proves that $E_M$ is in fact unital.
\epr

\btm\label{commutative-is-group}
Let $M\ne\{0\}$ be commutative and $M\simeq \hat N$ for some $N$. Then there is a locally compact group $G$ such that $M\simeq C_0(G)^{**}$.
\etm
\bpr
By construction, $M=A^{**}$ for the \cst-algebra $A=C^*(N_*^0)$, which is obviously commutative. Let $G$ be the spectrum of $A$, so that $A\simeq C_0(G)$. This is a locally compact space in the topology $\cal T_A$ induced by the Gelfand transform of $A$. Since $\Phi(N_*^0)$ is norm dense in $A$, the same topology is generated by the Gelfand transform of $\Phi(N_*^0)$.

Moreover, any subalgebra of $C(G)$ containing $C_0(G)$ generates the same topology $\cal T_A$ on $G$. Thus, the topology generated by $\Phi(N_*)$ is not stronger than $\cal T_A$, and with considerations above, it is equal to $\cal T_A$.

$G$ is identified with a subset in $N$, as the set of unitary characters of $N_*$ (by definition of $N_*^0$, they do not vanish on $N_*^0$). In fact, the characters of $N_*$ are given by group-like elements $u\in N$, i.e.\ such that $\Delta(u)=u\otimes u$, and a character is unitary if $u$ is unitary in $N$ and $\kappa(u)=u^*$. Then $G$ is a subgroup of the unitary group of $N$, so $G$ has a natural group structure. The group operations are continuous in the $N_*$-weak topology $\cal T_N$. But as we have shown above, it equals $\cal T_A$. Thus, the group operations are continuous on $G$ in  $\cal T_A$, so $G$ is a locally compact group.

Now, since $A\simeq C_0(G)$, we get $M=W^*(A)\simeq C_0(G)^{**}$. From Propositions \ref{dual-comultiplication} and \ref{dual-kappa} we conclude that this is also an isomorphism of coinvolutive Hopf--von Neumann algebras.
\epr

\bcor
If $\{0\}\ne M\simeq \hat{\hat M}$ and $M$ is commutative, then $M\simeq C_0(G)^{**}$ for a locally compact group $G$.
\ecor

\btm\label{cocommutative-is-group}
Let $M$ be cocommutative (i.e. $M_*$ is commutative) and $M\simeq \hat N$ for some $N$. Then, if $\Mhat\ne\{0\}$, there is a locally compact group $G$ such that $M\simeq W^*(G)$.
\etm
\bpr
As it is proved above, either $\Mhat=\{0\}$ or there is a locally compact group $G$ such that $\hat M\simeq C_0(G)^{**}$.
Assume further the second case. Then, by definition of the dual algebra, $G$ is the set of unitary characters of $M_*$ which do not vanish on $M_*^0$, so we can consider $G$ as a subset of $M$. More precisely, fix an imbedding $\i: G\to M$.

At the same time, every $g\in G$ is a character of $C_0(G)=C^*(M_*^0)$ and can be considered as an element of $\Mhat_*=M(G)$: this is, of course, the measure $\delta_g$ concentrated on $g$. In our notation, this means that $\hat\Phi_M(\delta_g)=\i(g)$. Under $\Phi_{\hat M}$, $\delta_g$ is mapped into $W^*(G)=\hat{\Mhat}$ in such a way that $D_M\Phi_{\hat M}(\delta_g)=\i(g)$ (by definition of $D_M$).

For the morphism $E_M: M\to \hat{\hat M}\simeq W^*(G)$, denote $u_g=E_M(\i (g))$. Since $\hat{\hat{\Mhat}}=C_0(G)^{**}\ne\{0\}$, $E_M$ is a unital *-homomorphism, so every $u_g$ is unitary. Since $E_M$ is a coalgebra morphism, we have $\Delta_{W^*(G)}(u_g)=(E_M\otimes E_M)(\Delta_M(\i (g)))=u_g\otimes u_g$, so $u_g$ is a character of $B(G)$. This together implies that $u_g=\Phi_{\hat M}(\delta_h)$ for some $h\in G$.

By Proposition \ref{DE=id}, $D_M(u_g)=D_M(E_M(\i(g))) = \i(g)$, but at the same time $D_M(u_g)=D_M(\Phi_{\hat M}(\delta_h))=\i(h)$. Thus, $g=h$. This implies that $E_M$ is in fact surjective, so this is an isomorphism, and the theorem is proved.
\epr

\section{\cst-algebraic version}\label{dual-cstar}

In this section we give a \cst-version of our construction. The class of \cst-algebras we use is different from coinvolutive Hopf \cst-algebras of Enock and Vallin \cite{enock-vallin} and from Hopf \cst-algebras of Vaes and Van Daele \cite{vaes-vandaele}. We will use a different term for this reason. Some ideas of \cite{vaes-vandaele} are used.

Here we will need definitions related to multiplier algebras, see Subsection \ref{section-multipliers}. Recall that for a homomorphism $\phi: A\to M(B)$, we denote by $\bar\phi$ its unique $(A^*, B^*)$--continuous extension to $M(A)$.

In the following definition the comultiplication is (as usual) a map $\Delta:A\to M(A\bar\otimes A)$. We can consider $M(A\bar\otimes A)$ as a subset of $W^*(A\bar\otimes A)$, and in this way we get a lifting of $\Delta$ to a map from $W^*(A)$ to $W^*(A\bar\otimes A)$. Composed with the canonical epimorphism $p:W^*(A\bar\otimes A)\to W^*(A)\bar\otimes W^*(A)$, this gives a morphism $\tilde\Delta:W^*(A)\to W^*(A)\bar\otimes W^*(A)$.

\begin{definition}\label{def-hopf-cstar}
A coinvolutive \cst-bialgebra is a \cst-algebra $A$ with a comultiplication $\Delta: A\to M(A\bar\otimes A)$ and coinvolution $\kappa:A\to A$ such that:\\
(i) $\Delta$ is a non-degenerate *-homomorphism, such that its lifting $\tilde\Delta: W^*(A)\to W^*(A)\bar\otimes W^*(A)$ is injective; \\
(ii) $\kappa$ is a *-antihomomorphism such that $\kappa^2={\rm id}$;\\
(iii) $\overline{(\Delta\otimes\id)}\Delta=\overline{(\id\otimes\Delta)}\Delta$ and $\overline{(\kappa\otimes\kappa)}=\bar\theta\Delta\kappa$, where $\theta$ is the flip map.
\end{definition}

It is well known that with the natural morphisms of \cst-algebras, the correspondence $G\mapsto C_0(G)$ for a locally compact group $G$ is not functorial: not every group homomorphism $\phi:G\to H$ gives by the pullback a morphism from $C_0(H)$ to $C_0(G)$. This justifies another definition: a morphism $\phi\in Mor(A,B)$ between two \cst-algebras is a homomorphism $\phi: A\to M(B)$. Some authors \cite{woron-mn} require even more: that $\phi$ is non-degenerate. But in our setting, as well as for morphisms of von Neumann algebras, it is better to allow degenerate morphisms too. This should agree with the Hopf structure as follows:

\begin{definition}\label{def-cstar-morphism}
Let $A$, $B$ be coinvolutive \cst-bialgebras. A morphism $\phi\in Mor(A,B)$ is a *-homomorphism $\phi:A\to M(B)$ such that\\
(i) $\bar\Delta_B \phi = (\overline{\phi\otimes\phi})\Delta_A$, where the range of $\overline{\phi\otimes\phi}$ is understood to be in $M(B)\bar\otimes M(B)\subset M(B\bar\otimes B)$;\\
(ii) $\bar\kappa_B \phi = \phi\kappa_A$.
\end{definition}

\bprop
Let $A$ be a coinvolutive \cst-bialgebra. Then the enveloping von Neumann algebra $W^*(A)$ of $A$ has a canonical structure of a coinvolutive Hopf--von Neumann algebra. If $\phi:A\to B$ is a morphism of coinvolutive \cst-bialgebras, then its normal lifting $\bar\phi: W^*(A)\to W^*(B)$ is a morphism of coinvolutive Hopf--von Neumann algebras.
\eprop
\bpr
 Since $\Delta$ is non-degenerate and $p$ is unital, $\tilde\Delta$ is also unital. Coassociativity of $\Delta$ together with uniqueness of these liftings implies that $\tilde\Delta$ is coassociative too.

Again by universality, $\kappa$ has a lifting to an anti-automorphism of $W^*(A)$. By continuity, \ref{def-hopf-cstar}(ii) implies the corresponding equality for $W^*(A)$.

It follows immediately from the identities \ref{def-cstar-morphism}(i,ii) that $\bar\phi$ is a morphism of coinvolutive Hopf--von Neumann algebras.
\epr

Next proposition is a simple fact on \cst-algebras. In the duality context it was applied, for example, in \cite{vaes-vandaele}.

\bprop\label{weak-and-norm}
Let $A$ be a \cst-algebra, and let $p_a\ge0$ be an increasing net in $A$. If $p_\a$ converges $A^*$-weakly in $A$, then it converges in norm.
\eprop
\bpr
Consider $p_\a$ as functions on the state space $S$ of $A$. Then they are continuous, nonnegative and increasing. Moreover, this net converges pointwise to some $p$. Since $S$ is compact, $p_\a\to p$ uniformly on it. This implies that in $A$, it converges in norm.
\epr

\btm
Let $A$ be a coinvolutive \cst-bialgebra and let $M$ be its enveloping (Hopf)--von Neumann algebra. Then $\hat A=C^*(M_*^0)$ has a canonical structure of a coinvolutive \cst-bialgebra.
\etm
\bpr
It is enough to prove the theorem in the nontrivial case $\hat A\ne\{0\}$.
We have seen in Section \ref{dual-von-neumann} that $W^*(\hat A)$ is a coinvolutive Hopf--von Neumann algebra, and that $\hat\kappa$ maps $\hat A$ to itself. We need next to prove that $\hat\Delta(\hat A)\subset M(\hat A\bar\otimes \hat A)$.

Recall that there is a map $\Phi: M_*\to \hat M$ extending the canonical inclusion of $M_*^0$ into its envelope. Suppose that $M$ acts on a Hilbert space $H$, and let $(e_\a)$ be an orthonormal basis of $H$. Denote $\phi_{\a\beta}=\Phi(\mu_{e_\beta e_\a})$.

Consider $\hat\Phi:\hat M_*\to M$ as a representation of $\hat M_*$. It is unitary by Proposition \ref{Phi-is-standard}. The coefficients of $\hat\Phi$ are $(\phi_{\a\beta})$. It follows that for all $\a$, we have $\hat A^*$-weak convergence in $\hat M$:
\beq\label{sum-phi}
\sum_{\beta} \phi_{\a\beta}\phi_{\a\beta}^*=\sum_\beta \phi_{\beta\a}^*\phi_{\beta\a}=1.
\eeq
By the \cite[1.4.2]{enock}, $\hat\Delta(\phi_{\a\a}) = \sum_{\beta} \phi_{\a\beta}\otimes \phi_{\beta\a}$ for every $\a$.

To show that $\hat\Delta(\phi_{\a\a})\in M(\hat A\bar\otimes \hat A)$, it suffices to show that $(a\otimes b)\hat\Delta(\phi_{\a\a})$ and $\hat\Delta(\phi_{\a\a})(a\otimes b)$ are both in $\hat A\bar\otimes \hat A$ for all $a,b\in \hat A$.

Fix $\a$ and denote $p_\beta=\phi_{\a\beta}$. Since $M_*^0$ is an ideal in $M_*$, for every $a\in \Phi(M_*^0)$ we have $ap_\beta\in \Phi(M_*^0)\subset \hat A$. By continuity, $\hat Ap_\beta\subset \hat A$, and similarly $p_\beta\hat A\subset \hat A$. Thus, $p_\beta\in M(\hat A)$.

For $a\in \hat A$, from \eqref{sum-phi} we conclude that $\sum_\beta a p_\beta p_\beta^*a^*=aa^*$ for every $a\in \hat A$, weakly in $\hat A$. Applying Proposition \ref{weak-and-norm}, we see that this series converges also in norm.

For $b\in \hat A$, we have $\phi_{\beta\a}^*b^*b\phi_{\beta\a}\le \|b\|^2 \phi_{\beta\a}^*\phi_{\beta\a}$.
 Then $\sum_\beta \phi_{\beta\a}^*b^*b\phi_{\beta\a}\le \|b\|^2\sum_\beta \phi_{\beta\a}^*\phi_{\beta\a} = \|b\|^2$, so the net of partial sums of this series is bounded.

Now we see that the series $\sum_{\beta} a\phi_{\a\beta}\otimes b\phi_{\beta\a}$ converges in the Haagerup norm on $\hat A\otimes \hat A$: for every finite subset $B$ of indices,
$$
\|\sum_{\beta\in B} a\phi_{\a\beta}\otimes b\phi_{\beta\a}\|_h^2 \le \|\sum_{\beta\in B} a\phi_{\a\beta}\phi_{\a\beta}^*a^*\|\,\|\sum_{\gamma\in B} \phi_{\gamma\a}^*b^*b\phi_{\gamma\a}\|.
$$
Hence, $(a\otimes b)\hat\Delta(\phi_{\a\a})\in \hat A\otimes_h\hat A\subset \hat A\bar\otimes \hat A$. The other inclusion is proved identically.
Thus, $\hat\Delta(\phi_{\a\a})$ is in $M(\hat A\bar\otimes \hat A)$ (and even in $M(\hat A\otimes_h\hat A)$). By polarization, we get this inclusion also for $\hat\Delta(\phi_{\a\beta})$ for all $\a,\beta$. Since the set of all $\mu_{e_\beta e_\a}$ is total in $M_*$ and as a consequence the set of all $\phi_{\a\beta}$ is weakly total in $\Mhat$, we get by continuity that $\hat\Delta(\hat A)\subset \hat\Delta(\Mhat)$ is also in $M(\hat A\bar\otimes \hat A)$. Clearly $\hat\Delta$ is non-degenerate.

The equalities \ref{def-hopf-cstar}(ii) we have for granted, since they are valid in $W^*(A)$.
\epr

\begin{definition}
Let $A$ be a coinvolutive \cst-bialgebra. The coinvolutive \cst-bialgebra $\hat A=C^*(M_*^0)$ will be called the dual coinvolutive \cst-bialgebra of $A$.
\end{definition}

Now we can translate the von Neumann algebraic duality theorem into the language of \cst-algebras.

\btm
Let $\cal {CCB}_0$ be the category of coinvolutive \cst-bialgebras, and let $\cal {CCB}$ be the full subcategory of $\cal {CCB}_0$ which has as its objects all $A\in \cal {CCB}_0$ such that $A\simeq \hat{\hat A}$. Then $\hat{\;}$ is a contravariant functor on $\cal {CCB}$, such that:\\
\indent$\hat{C_0(G)}=C^*(G)$ and $\hat{C^*(G)}=C_0(G)$, for every locally compact group $G$;\\
\indent$\hat A\simeq (\hat A)\widehat{\widehat\;}$ for every Kac \cst-algebra $A$;\\
\indent if $A\ne\{0\}$ is commutative and $A\simeq \hat B$ for some $B$, then $A$ is isomorphic to $C_0(G)$ for some locally compact group $G$;\\
\indent if $\hat A\ne\{0\}$, $A$ is co-commutative and $A\simeq \hat B$ for some $B$, then $A$ is isomorphic to $C^*(G)$ for some locally compact group $G$.
\etm
\bpr
The only statement to check is that the dual of a morphism is always well defined. Let $\phi:A\to M(B)$ be a morphism with $A,B\in \cal {CCB}$. Denote $\frak A=W^*(A)$, $\frak B=W^*(B)$. Composing with the canonical imbedding $i:M(B)\to \frak B$, we get a homomorphism $i\circ\phi:A\to \frak B$ and its normal lifting $\bar\phi:\frak A\to \frak B$.

Consider the dual (in the sense of Proposition \ref{dual-morphism}) morphism $\hat\phi: \hat{\frak B}\to \hat{\frak A}$.
By definition, $\hat\phi\circ\Phi_{\frak B} = \Phi_{\frak A}\circ\bar\phi_*$ on $\frak B_*^0$. It follows that $\hat\phi(b)\in \Phi_{\frak A}(\frak A_*)$ for every $b\in \Phi_{\frak B}(\frak B_*^0)$. Now, $\Phi_{\frak A}(\a)\hat A\subset \hat A$ for every $\a\in \frak A_*$, since $\frak A_*^0$ is an ideal in $\frak A_*$. This implies that $\hat\phi(\hat B)$ is in $M(\hat A)$, and the rest is obvious.
\epr

In conclusion, we will draw one more example of a dual algebra: $A=SU_2(n)$, considered with the unitary antipode. The dual is in this case very degenerate, $\hat A=\C$. If we used the usual antipode, the second dual would be again $A$, but this would require changing the definitions as to allow the unbounded antipode. This will be done elsewhere.

\bex
Let $A=SU_q(2)$ be the Woronowicz's quantum $SU(2)$ group \cite{woron-suq}, $q\in\R$, $q\ne0$. As every compact quantum group, $A$ is the closed linear span of the coefficients $u^\a_{ij}$ of irreducible finite-dimensional corepresentations $\pi^\a$ of $A$. In the case of $SU_q(2)$, one can index them by $\a\in \N_0$, so that the dimension of $\pi^\a$ is $2\a+1$ (see, e.g., \cite[4.2]{klimyk}).

The usual antipode $S$ is given by $S(u^\a_{ij}) =  u^{\a*}_{ji}$, and does not extend to a bounded map on $A$. However, one can define \cite[\S~5]{mn} a unitary map $\kappa$ on $A$ by $\kappa(u^\a_{ij}) = q^{j-i}u^{\a*}_{ji}$ 
such that $A$ becomes a coinvolutive \cst-bialgebra in the sense above (see \cite[4.2.4]{klimyk} to extend the formula for $\kappa$ from the fundamental representation to any $\a$).

Let $h$ be the Haar state on $A$. In $A^*$, we have a family of coordinate functionals $e^\a_{ij}$ such that $e^\a_{ij}(u^\beta_{kl})=\delta_{\a\beta}\delta_{ik}\delta_{jl}$. As it follows from the orthogonality of $u^\a_{ij}$ with respect to the scalar product $\langle a,b\rangle=h(b^*a)$ \cite[Theorem 17]{klimyk}, they may be expressed in terms of $h$ as
\beq\label{e-via-h}
e^\a_{ij}(x)=h((u^\a_{ij})^* x)/h((u^\a_{ij})^*u^\a_{ij}).
\eeq
These functionals separate points of $A$: for $x\in A$, vanishing on every $e^\a_{ij}$ means being orthogonal to every $u^\a_{ij}$, and it is known that the set of $u^\a_{ij}$ is total in $A$ with $\langle \cdot,\cdot\rangle$.
Denote by $I_\a$ the linear span of $e^\a_{ij}$ for fixed $\a$. This is a finite-dimensional ideal in $A^*$, so every irreducible representation $\pi$ is either zero or irreducible on $I_\a$. In the latter case it quickly follows that $\pi=\pi^\a$. Moreover, if a representation vanishes on every $I_\a$ then it is zero since $\cup I_\a$ separates points of $A$.

If we calculate the involution of $e^\a_{ij}$ using the antipode $\kappa$, we get that $(e^\a_{ij})^*=q^{i-j}e^\a_{ji}$. At the same time, $\pi^\a(e^\a_{ij})^*=\pi^\a(e^\a_{ji})$, so that $\pi^\a$ is not involutive if $\dim\pi^\a>1$. Thus, $\pi^\a$ is unitary only if $\pi^\a(\mu)=\mu(1)$.

Thus, we have only one (one-dimensional) unitary representation, and non-unitary representations do not separate points of $A^*$; thus, $\hat A=\C$.
\eex


\begin{thebibliography}{99}

\bibitem{baaj} Baaj S., Skandalis G. Unitaires multiplicatifs et dualit\'e pour les produits crois\'es de \cst-alg\`ebres. {\it Ann. Sci. \'Ecole Norm. Sup. (4)} {\bf 26} (1993), no.~4, 425--488.
\bibitem{blecher} Blecher D. P., Le Merdy C. {\it Operator algebras and their modules: an operator space approach}. Oxford University Press, 2004.
\bibitem{dixmier} Dixmier J. {\it \cst-algebras}. North Holland, 1977.
\bibitem{effros} Effros E.~G., Ruan Zh.-J. {\it Operator spaces}. Oxford, 2000.
\bibitem{enock} Enock M., Schwartz J.-M. {\it Kac algebras and duality of locally compact groups}. Springer, 1992.
\bibitem{enock-vallin} Enock M., Vallin J.-M. \cst-alg\`ebres de Kac et alg\`ebres de Kac. {\it Proc. London Math. Soc. (3)} {\bf 66} (1993), 619--650.
\bibitem{hajac} {\it Quantum symmetry in noncommutative geometry}. Hajac P. M. (Ed.) EMS, to appear.
\bibitem{HR} Hewitt E., Ross~K.~A. {\it Abstract harmonic analysis I, II.} Springer, 3rd printing, 1997.
\bibitem{kadison} Kadison R. V., Ringrose J. R. {\it Fundamentals of the theory of operator algebras. Vol. I, II}. Academic Press, 1986.
\bibitem{kaniuth} Kaniuth E. {\it A course in commutative Banach algebras}. Springer, 2009.
\bibitem{kirchberg} Kirchberg E. {\it Darstellungen coinvolutiver Hopf-$W^*$-Algebren und ihre Anwendung in der nicht-abelschen Dualit\"atstheorie lokalkompakter Gruppen}. Thesis, Berlin, 1977.
\bibitem{klimyk} Klimyk A., Schm\"udgen K. {\it Quantum groups and their representations}. Springer, 1997.
\bibitem{kustermans} J. Kustermans, Locally compact quantum groups in the universal setting. {\it Internat. J. Math.} {\bf 12} (2001), no. 3, 289--338.
\bibitem{woron-mn} Masuda T., Nakagami Y., Woronowicz S. L. A \cst-algebraic framework for quantum groups. {\it Internat. J. Math.} {\bf 14} (2003), no.~9, 903--1001.
\bibitem{mn} Masuda T., Nakagami Y. A von Neumann algebra framework for the duality of the quantum groups. {\it Publ. RIMS, Kyoto Univ.} {\bf 30} (1994), 799--850.
\bibitem{ng} Ng Ch.-K. Duality of Hopf \cst-algebras. {\it Int. J. Math.} {\bf 13} no.~9 (2002), 1009--1025.
\bibitem{palmer} Palmer Th. W. {\it Banach algebras and the general theory of *-algebras. Vol. II}. Cambridge University Press, 2001.
\bibitem{pedersen} Pedersen G. K. {\it \cst-algebras and their automorphism groups}, Academic Press (1979).
\bibitem{takesaki-tohoku-1958} Takesaki M. On the conjugate space of operator algebra. {\it Tohoku Math. J.}, {\bf 10} (1958), 194--203.
\bibitem{taylor} Taylor J. L. $L$-Subalgebras of $M(G)$. {\it Trans. AMS} {\bf 135} (1969), 105--113.
\bibitem{vaes-vandaele} Vaes S., Van Daele A. Hopf C*-Algebras. {\it Proc. London Math. Soc.} {\bf 82} (2): 337--384 (2001).
\bibitem{walter} Walter M.~E. $W^*$-algebras and nonabelian harmonic analysis. {\it J. Funct. Anal.} {\bf 11} (1972), 17--38.
\bibitem{woron} S.~L.~Woronowicz. From multiplicative unitaries to quantum groups. {\it Internat. J. Math.} {\bf 7} (1996), no.~1, 127--149.
\bibitem{woron-suq} S.~L.~Woronowicz. Tannaka-Krein duality for compact matrix pseudogroups. Twisted $SU(N)$ groups. {\it Invent. Math.} {\bf 93} (1988), 35--76.
\end{thebibliography}
\end{document}